\documentclass[10pt,reqno]{amsart}
\usepackage{latexsym,amsmath,amssymb,amscd,amsfonts,amsthm,mathrsfs}
\usepackage[all]{xy}

\def\today{\ifcase \month \or
   January \or February \or March \or April \or
   May \or June \or July \or August \or
   September \or October \or November \or December \fi
   \space\number\day , \number\year}
   
\makeatletter
  \newcommand\@dotsep{4.5}
  \def\@tocline#1#2#3#4#5#6#7{\relax
     \ifnum #1>\c@tocdepth 
     \else
     \par \addpenalty\@secpenalty\addvspace{#2}%
     \begingroup \hyphenpenalty\@M
     \@ifempty{#4}{%
     \@tempdima\csname r@tocindent\number#1\endcsname\relax
        }{%
         \@tempdima#4\relax
           }%
      \parindent\z@ \leftskip#3\relax \advance\leftskip\@tempdima\relax
      \rightskip\@pnumwidth plus1em \parfillskip-\@pnumwidth
       #5\leavevmode\hskip-\@tempdima #6\relax
       \leaders\hbox{$\m@th
       \mkern \@dotsep mu\hbox{.}\mkern \@dotsep mu$}\hfill
       \hbox to\@pnumwidth{\@tocpagenum{#7}}\par
       \nobreak
        \endgroup
         \fi}
\makeatother 

\begin{document}
\parskip=5pt


\makeatletter
\@addtoreset{figure}{section}
\def\thefigure{\thesection.\@arabic\c@figure}
\def\fps@figure{h,t}
\@addtoreset{table}{bsection}

\def\thetable{\thesection.\@arabic\c@table}
\def\fps@table{h, t}
\@addtoreset{equation}{section}
\def\theequation{
\arabic{equation}}
\makeatother

\newcommand{\bfi}{\bfseries\itshape}

\newtheorem{theorem}{Theorem}
\newtheorem{acknowledgment}[theorem]{Acknowledgment}
\newtheorem{corollary}[theorem]{Corollary}
\newtheorem{definition}[theorem]{Definition}
\newtheorem{example}[theorem]{Example}
\newtheorem{lemma}[theorem]{Lemma}
\newtheorem{notation}[theorem]{Notation}
\newtheorem{proposition}[theorem]{Proposition}
\newtheorem{remark}[theorem]{Remark}
\newtheorem{setting}[theorem]{Setting}

\numberwithin{theorem}{section}
\numberwithin{equation}{section}

\renewcommand{\1}{{\bf 1}}
\newcommand{\Ad}{{\rm Ad}}
\newcommand{\Alg}{{\rm Alg}\,}
\newcommand{\Aut}{{\rm Aut}\,}
\newcommand{\ad}{{\rm ad}}
\newcommand{\Borel}{{\rm Borel}}
\newcommand{\card}{{\rm card}\,}
\newcommand{\Ci}{{\mathscr C}^\infty}
\newcommand{\Cpol}{{\mathscr C}^\infty_{\rm pol}}
\newcommand{\Der}{{\rm Der}\,}
\newcommand{\de}{{\rm d}}
\newcommand{\ee}{{\rm e}}
\newcommand{\End}{{\rm End}\,}
\newcommand{\ev}{{\rm ev}}
\newcommand{\hotimes}{\widehat{\otimes}}
\newcommand{\id}{{\rm id}}
\newcommand{\ie}{{\rm i}}
\newcommand{\GL}{{\rm GL}}
\newcommand{\gl}{{{\mathfrak g}{\mathfrak l}}}
\newcommand{\Hom}{{\rm Hom}\,}
\newcommand{\Img}{{\rm Im}\,}
\newcommand{\Ind}{{\rm Ind}}
\newcommand{\Ker}{{\rm Ker}\,}
\newcommand{\Lie}{\text{\bf L}}
\newcommand{\local}{{\rm loc}}
\newcommand{\m}{\text{\bf m}}
\newcommand{\pr}{{\rm pr}}
\newcommand{\Ran}{{\rm Ran}\,}
\renewcommand{\Re}{{\rm Re}\,}
\newcommand{\so}{\text{so}}
\newcommand{\spa}{{\rm span}\,}
\newcommand{\supp}{{\rm supp}\,}
\newcommand{\Tr}{{\rm Tr}\,}
\newcommand{\Op}{{\rm Op}}
\newcommand{\U}{{\rm U}}

\newcommand{\CC}{{\mathbb C}}
\newcommand{\RR}{{\mathbb R}}
\newcommand{\TT}{{\mathbb T}}

\newcommand{\Ac}{{\mathscr A}}
\newcommand{\Bc}{{\mathscr B}}
\newcommand{\Cc}{{\mathscr C}}
\newcommand{\Dc}{{\mathscr D}}
\newcommand{\Ec}{{\mathscr E}}
\newcommand{\Fc}{{\mathscr F}}
\newcommand{\Hc}{{\mathscr H}}
\newcommand{\Jc}{{\mathscr J}}
\renewcommand{\Mc}{{\mathscr M}}
\newcommand{\Nc}{{\mathscr N}}
\newcommand{\Oc}{{\mathscr O}}
\newcommand{\Pc}{{\mathscr P}}
\newcommand{\Rc}{{\mathscr R}}
\newcommand{\Sc}{{\mathscr S}}
\newcommand{\Tc}{{\mathscr T}}
\newcommand{\Vc}{{\mathscr V}}
\newcommand{\Uc}{{\mathscr U}}
\newcommand{\Xc}{{\mathscr X}}
\newcommand{\Yc}{{\mathscr Y}}
\newcommand{\Wig}{{\mathscr W}}

\newcommand{\Bg}{{\mathfrak B}}
\newcommand{\Fg}{{\mathfrak F}}
\newcommand{\Gg}{{\mathfrak G}}
\newcommand{\Ig}{{\mathfrak I}}
\newcommand{\Jg}{{\mathfrak J}}
\newcommand{\Lg}{{\mathfrak L}}
\newcommand{\Pg}{{\mathfrak P}}
\newcommand{\Sg}{{\mathfrak S}}
\newcommand{\Xg}{{\mathfrak X}}
\newcommand{\Yg}{{\mathfrak Y}}
\newcommand{\Zg}{{\mathfrak Z}}

\newcommand{\ag}{{\mathfrak a}}
\newcommand{\bg}{{\mathfrak b}}
\newcommand{\dg}{{\mathfrak d}}
\renewcommand{\gg}{{\mathfrak g}}
\newcommand{\hg}{{\mathfrak h}}
\newcommand{\kg}{{\mathfrak k}}
\newcommand{\mg}{{\mathfrak m}}
\newcommand{\n}{{\mathfrak n}}
\newcommand{\og}{{\mathfrak o}}
\newcommand{\pg}{{\mathfrak p}}
\newcommand{\sg}{{\mathfrak s}}
\newcommand{\tg}{{\mathfrak t}}
\newcommand{\ug}{{\mathfrak u}}
\newcommand{\zg}{{\mathfrak z}}

\newcommand{\ZZ}{\mathbb Z}
\newcommand{\NN}{\mathbb N}
\newcommand{\BB}{\mathbb B}

\newcommand{\ep}{\varepsilon}

\newcommand{\hake}[1]{\langle #1 \rangle }

\newcommand{\scalar}[2]{\langle #1 ,#2 \rangle }
\newcommand{\vect}[2]{(#1_1 ,\ldots ,#1_{#2})}
\newcommand{\norm}[1]{\Vert #1 \Vert }
\newcommand{\normrum}[2]{{\norm {#1}}_{#2}}

\newcommand{\upp}[1]{^{(#1)}}
\newcommand{\p}{\partial}

\newcommand{\opn}{\operatorname}
\newcommand{\slim}{\operatornamewithlimits{s-lim\,}}
\newcommand{\sgn}{\operatorname{sgn}}

\newcommand{\seq}[2]{#1_1 ,\dots ,#1_{#2} }
\newcommand{\loc}{_{\opn{loc}}}

\makeatletter
\title[Smooth vectors and Weyl-Pedersen calculus]{Smooth vectors 
and Weyl-Pedersen calculus for 
representations of nilpotent Lie groups}
\author{Ingrid Belti\c t\u a 
and Daniel Belti\c t\u a
}
\address{Institute of Mathematics ``Simion Stoilow'' 
of the Romanian Academy, 
P.O. Box 1-764, Bucharest, Romania}
\email{Ingrid.Beltita@imar.ro}
\email{Daniel.Beltita@imar.ro}
\keywords{Weyl calculus; nilpotent Lie group; semidirect product}
\subjclass[2000]{Primary 22E25; Secondary 22E27, 35S05, 47G30}
\dedicatory{To Professor Ion Colojoar\u a for his eightieth birthday}
\date{October 25, 2009}
\makeatother

\begin{abstract} 
We present some recent results on smooth vectors for 
unitary irreducible representations of nilpotent Lie groups. 
Applications to the Weyl-Pedersen calculus of pseudo-differential operators 
with symbols on the coadjoint orbits are also discussed.  
\end{abstract}

\maketitle

\section{Introduction}\label{Sect1}

`Weyl-Pedersen calculus' is the name proposed in \cite{BB09c} 
for the remarkable correspondence $a\mapsto\Op^\pi(a)$ constructed 
by N.V.~Pedersen in \cite{Pe94} as a generalization of 
the pseudo-differential Weyl calculus on~$\RR^n$. 
Here $\pi\colon G\to\Bc(\Hc)$ is any unitary irreducible representation 
of a connected, simply connected, nilpotent Lie group $G$, 
the symbol $a$ can be any tempered distribution on 
the coadjoint orbit~$\Oc$ corresponding to $\pi$ by the orbit method 
of~\cite{Ki62}, and $\Op^\pi(a)$ is a linear operator 
in the representation space $\Hc$, which is in general unbounded. 
When $\pi$ is the Schr\"odinger representation of 
the $(2n+1)$-dimensional Heisenberg group, 
the correspondence $a\mapsto\Op^\pi(a)$ is precisely 
the calculus suggested by H.~Weyl in \cite{We28} 
for applications to quantum mechanics. 
This calculus was later developed by L.~H\"ormander in \cite{Hor79} 
and made into a powerful calculus of pseudo-differential operators on~$\RR^n$; 
see \cite{Hor07}. 

In the present paper we discuss the classical notion of smooth vectors 
---and the related notion of smooth operators--- 
with a view toward their crucial importance for the Weyl-Pedersen calculus. 
We then approach a related circle of ideas that recently emerged in \cite{BB09c}, 
namely the modulation spaces 
for unitary irreducible representations of nilpotent Lie groups. 
We take the opportunity of this discussion 
to extend some known facts to a setting 
where they hold true in a natural degree of generality 
(see for instance Theorem~\ref{smooth_ideal} below). 
We also take a close look at some new examples of unitary irreducible representations 
and find out their related notions which illustrate the main theme 
of the present paper: the preduals of the corresponding coadjoint orbits, 
their ambiguity function, or their space of smooth vectors 
(see Proposition~\ref{double} and Corollary~\ref{double_cor}). 

Let us mention that the importance of the Weyl-Pedersen calculus and the related circle of ideas  
goes far beyond the framework of representation theory of nilpotent Lie groups. 
Many other interesting developments within 
the theory of partial differential equations and the finite-dimensional Lie theory 
can be found for instance in the references \cite{An69}, \cite{An72}, 
\cite{How80}, \cite{Mi82}, \cite{Me83}, \cite{How84}, \cite{HRW84}, \cite{HN85}, \cite{Mi86}, 
\cite{FG92}, \cite{Ma91}, \cite{Ma95}, 
and \cite{Ma07}. 
Moreover, one can use a similar construction 
even for representations of certain \emph{infinite-dimensional Lie groups} 
in order to provide a geometric explanation for 
the gauge covariance for the magnetic Weyl calculus 
of \cite{MP04}, \cite{IMP07}, \cite{IMP09}, \cite{MP09} and the references therein. 
The representation theoretic approach to the magnetic Weyl calculus 
has been taken up in the papers \cite{BB09a}, \cite{BB09b}; 
see also the survey \cite{BB09d}.

The structure of the present paper is summarized in the following table of contents: 

\tableofcontents

\noindent\textbf{Notation and background.}
Throughout the paper we denote by $\Sc(\Vc)$ the Schwartz space 
on a finite-dimensional real vector space~$\Vc$. 
That is, $\Sc(\Vc)$ is the set of all smooth functions 
that decay faster than any polynomial together with 
their partial derivatives of arbitrary order. 
Its topological dual ---the space of tempered distributions on $\Vc$--- 
is denoted by $\Sc'(\Vc)$. 
We shall also have the occasion to 
use these notions with $\Vc$ replaced by a coadjoint orbit 
of a nilpotent Lie group. 
In this situation we need the notion of polynomial structure 
on a manifold; see Sect.~1 in \cite{Pe89} for details. 
We use $\scalar{\cdot}{\cdot}$ to denote any duality pairing between 
finite-dimensional real vector spaces whose meaning is clear 
from the context. 

We shall also use the convention that the Lie groups are denoted by 
upper case Latin letters and the Lie algebras are denoted 
by the corresponding lower case Gothic letters.

As regards the background information for the present paper, 
we refer to \cite{Hor07}, \cite{Fo89}, and \cite{Gr01} for 
basic notions of pseudo-differential Weyl calculus on $\RR^n$. 
The necessary notions of representation theory for nilpotent Lie groups 
(in particular, 
the correspondence between the coadjoint orbits and the unitary irreducible representations)
can be found in \cite{Pu67}, \cite{CG90}, and \cite{Ki04}; 
see also \cite{Wa72} and \cite{Ki76}. 
Our references for topological vector spaces, nuclear spaces, and related topics 
are \cite{Sch66}, \cite{Tr67}, and \cite{Co68}. 

\section{Smooth vectors for Lie group representations}

The smooth vectors have been a basic tool in representation theory 
of Lie groups; see for instance the early paper \cite{Ga47} and 
the classical monographs \cite{Wa72} and \cite{Ki76}. 
In this section we discuss some of the very basic properties 
of the smooth vectors for the purpose of providing the necessary 
background for the later developments in the present paper. 

\begin{notation}
\normalfont
Throughout this section we shall use the following notation:
\begin{itemize}
\item $G$ is a connected unimodular Lie group with the Lie algebra $\gg$; 
\item $\de x$ denotes a fixed Haar measure on $G$; 
\item $\Vc$ and $\Yc$ are some complex Banach spaces; 
\item $\pi\colon G\to \Bc(\Yc)$ is a representation 
which is continuous, in the sense that for every $x\in \Yc$ 
the mapping $\pi(\cdot)x\colon G\to \Yc$ is continuous. 
\end{itemize}
\qed
\end{notation}

\subsection{Distribution theory on Lie groups} 
Some references for distribution theory on Lie groups 
are \cite{Br56} and \cite{Wa72}. 
The present subsection just records a few basic notions and properties 
needed later. 

\begin{definition}
\normalfont
We define the \emph{spaces of test functions} on the Lie group $G$ as follows: 
\begin{enumerate}
\item The space 
$$\Ec(G,\Vc):=\{\phi\colon G\to\Vc\mid \phi\text{ is smooth}\}$$
with the usual topology of a Fr\'echet space 
(given by the uniform convergence 
on compact sets of functions and their partial derivatives in local charts).
\item The space 
$$\Dc(G,\Vc):=\{\phi\in\Ec(G,\Vc)\mid \supp\,\phi\text{ is compact}\}$$
with the usual topology of an inductive limit of Fr\'echet spaces. 
\end{enumerate}
If $\Vc=\CC$ then we denote simply $\Ec(G,\CC)=\Ec(G)$ and $\Dc(G,\CC)=\Dc(G)$. 
For every integer $m\ge1$ we shall also need the function space 
$$\Cc^m_0(G):=\{\phi\colon G\to\CC\mid \phi\text{ is of class $\Cc^m$ and }\supp\,\phi\text{ is compact}\} $$
with its usual topology of an inductive limit of Banach spaces. 

We then define the \emph{spaces of vector valued distributions}
$$\Dc'^{\Yc}(G,\Vc):=\{u\colon\Dc(G,\Vc)\to\Yc\mid u\text{ is linear and continuous}\}$$
and 
$$\Ec'^{\Yc}(G,\Vc):=\{u\colon\Ec(G,\Vc)\to\Yc\mid u\text{ is linear and continuous}\}$$
and endow them with the topology of pointwise convergence. 
We denote the evaluation mapping by 
$$\langle\cdot,\cdot\rangle\colon
\Dc'^{\Yc}(G,\Vc)\times\Dc(G,\Vc)\to\Yc,\quad 
\langle u,\phi\rangle:=u(\phi), $$
and similarly for 
$\langle\cdot,\cdot\rangle\colon
\Ec'^{\Yc}(G,\Vc)\times\Ec(G,\Vc)\to\Yc$. 

For $\Vc=\CC$ we denote simply
$\Ec'^{\Yc}(G,\CC)=\Ec'^{\Yc}(G)$ and $\Dc'^{\Yc}(G,\CC)=\Dc'^{\Yc}(G)$. 
If also $\Yc=\CC$, then we further denote 
$\Ec'^{\Yc}(G)=\Ec'(G)$ and $\Dc'^{\Yc}(G)=\Dc'(G)$.
\qed
\end{definition}

\begin{definition}
\normalfont
The \emph{support} of the distribution $u\in\Dc'^{\Yc}(G,\Vc)$ 
is denoted by $\supp u$ and is defined as the intersection 
of all the closed sets $F\subseteq G$ such that for every $\phi\in\Dc(G,\Vc)$ 
with $F\cap\supp\phi=\emptyset$ we have $\langle u,\phi\rangle=0$. 
\qed
\end{definition}

\begin{remark}\label{funct}
\normalfont
Let $L^1_{\loc}(G)$ denote the linear space of (equivalence classes of) 
measurable functions on $G$ which are absolutely integrable with respect to the Haar measure $\de x$ on every compact subset of~$G$. 
Then there exists a natural linear embedding 
$L^1_{\loc}(G)\hookrightarrow\Dc'(G)$. 
Specifically, every $f\in L^1_{\local}(G)$ gives rise to 
a distribution also denoted by $f$ and defined by 
$$(\forall \phi\in\Dc(G))\quad \langle f,\phi\rangle=\int\limits_G f\phi\de x. $$
Note that $L^1_{\local}(G)$ contains many function spaces on $G$, 
like the space of continuous functions, 
or the space of smooth functions $\Ec(G)$, 
or the Lebesgue space $L^p(G)$ if $1\le p\le\infty$ etc.
\qed
\end{remark}

\begin{remark}\label{comp_supp}
\normalfont
We have 
$$\Ec'(G)=\{u\in\Dc'(G)\mid\supp u\text{ is compact}\}. $$
For every compact set $K\subseteq G$ we denote 
$\Ec'_K(G)=\{u\in\Dc'(G)\mid\supp u\subseteq K\}. $
\qed
\end{remark}

\begin{remark}\label{tens}
\normalfont
We recall from \cite{Sch66} and \cite{Tr67} that 
the locally convex spaces $\Ec(G)$ and $\Dc(G)$ 
are nuclear. 
Moreover, we have the linear topological isomorphisms 
$$\Ec(G,\Yc)\simeq\Ec(G)\hotimes\Yc
\text{ and }\Dc(G,\Yc)\simeq\Dc(G)\hotimes\Yc,$$
which are natural in the sense that every pair $(\phi,y)\in\Ec(G)\times\Yc$ 
corresponds to the function $\phi(\cdot)y\in\Ec(G,\Yc)$. 
Also recall the the linear topological isomorphisms 
$$\Ec(G)\hotimes\Ec(G)\simeq\Ec(G\times G) 
\quad\text{ and }\quad
\Dc(G)\hotimes\Dc(G)\simeq\Dc(G\times G) $$
that take a pair $(\phi_1,\phi_2)\in\Dc(G)\times\Dc(G)$ 
to the function $\phi_1\otimes\phi_2$ defined by $(x_1,x_2)\mapsto \phi_1(x_1)\phi_2(x_2)$. 
\qed
\end{remark}

\begin{example}
\normalfont 
Here are some examples of vector valued distributions 
that will be needed in the sequel. 
\begin{enumerate}
\item For arbitrary $g\in G$
the $\Yc$-valued \emph{Dirac distribution} $\delta_g^{\Yc}\in\Ec'^{\Yc}(G,\Yc)$ is defined by 
$$\delta_g^{\Yc}\colon\Ec(G,\Yc)\to\Yc, \quad 
\langle\delta_g^{\Yc},\phi\rangle=\phi(g).  $$
If $\Yc=\CC$ then we denote simply $\delta_g^{\Yc}=\delta_g$. 
\item By using Remark~\ref{tens}, one can define a canonical linear mapping 
$$\Ec'(G)\to\Ec'^{\Yc}(G,\Yc),\quad u\mapsto u\otimes\id_{\Yc}, $$
which takes every distribution $u\colon\Ec(G)\to\CC$ 
to its tensor product with the identity operator $\id_{\Yc}\colon\Yc\to\Yc$. 
\end{enumerate}
\qed
\end{example}

\begin{definition}\label{conv_def}
\normalfont 
Let $u_1,u_2\in\Ec'(G)$. 
Then the \emph{tensor product of distributions} 
$u_1\otimes u_2\in\Ec'(G\times G)$ 
can be defined by using Remark~\ref{tens} 
such that $\langle u_1\otimes u_2,\phi_1\otimes\phi_2\rangle
=\langle u_1,\phi_1\rangle\cdot\langle u_2,\phi_2\rangle$. 
On the other hand, there exists a continuous linear co-product 
$$\Ec(G)\to\Ec(G\times G),\quad \phi\mapsto\phi^\Delta, $$
where $\phi^\Delta(x_1,x_2):=\phi(x_1x_2)$ whenever $x_1,x_2\in G$ 
and $\phi\in\Ec(G)$. 
The \emph{convolution product of distributions} $u_1\ast u_2\in\Ec'(G)$ 
is defined by 
$$(\forall \phi\in\Ec(G))\quad 
\langle u_1\ast u_2,\phi\rangle:=\langle u_1\otimes u_2,\phi^\Delta\rangle.$$
The convolution product makes the distribution space $\Ec'(G)$ 
into an associative algebra 
whose unit element is the Dirac distribution $\delta_{\1}\in\Ec'(G)$. 
\qed
\end{definition}

\begin{example}\label{conv_ex}
\normalfont 
Let us consider a few simple properties of the convolution product. 
\begin{enumerate}
\item For arbitrary $g_1,g_2\in G$ we have $\delta_{g_1}\ast\delta_{g_2}=\delta_{g_1g_2}$. 
\item For every $u_1,u_2\in\Ec'(G)$ we have 
$$\supp(u_1\ast u_2)\subseteq\{x_1x_2\mid x_j\in\supp u_j\text{ for }j=1,2\}.$$
\end{enumerate}
\qed
\end{example}

\begin{definition}\label{univ_def}
\normalfont
We shall think of $\gg$ as a real subalgebra of its \emph{complexification} 
$\gg_{\CC}:=\CC\otimes_{\RR}\gg$, hence $\gg_{\CC}=\gg\dotplus\ie\gg$. 
The \emph{universal enveloping algebra} $\U(\gg_{\CC})$ is 
the  complex unital associative algebra satisfying the following conditions: 
\begin{enumerate}
\item\label{univ_def_item1} 
The complexification $\gg_{\CC}$ is a Lie subalgebra of $\U(\gg_{\CC})$. 
\item\label{univ_def_item2} 
For every complex unital associative algebra $\Ac$ and every linear mapping 
$\theta\colon\gg_{\CC}\to\Ac$ satisfying  $\theta([X,Y])=\theta(X)\theta(Y)-\theta(Y)\theta(X)$ 
for all $X,Y\in\gg_{\CC}$ there exists a unique extension of 
$\theta$ to a homomorphism of complex unital associative algebras $\U(\gg_{\CC})\to\Ac$. 
\end{enumerate}
One can prove that there always exists an algebra $\U(\gg_{\CC})$ satisfying 
these conditions and it is uniquely determined up to an isomorphism of 
complex unital associative algebras. 
Moreover, there exists a unique (complex-)linear mapping $\U(\gg_{\CC})\to\U(\gg_{\CC})$, $u\mapsto u^\perp$ 
such that 
$$(vw)^\perp=w^\perp v^\perp,\quad (w^\perp)^\perp=w,
\quad\text{and}\quad X^\perp=-X $$
for every $u,v\in\U(\gg_{\CC})$ and $X\in\gg$. 
(See \cite{Dix74} for more details). 
\qed
\end{definition}

\begin{example}\label{univ_ex}
\normalfont
If $\gg$ is an abelian Lie algebra of dimension $n$, 
then $\U(\gg_{\CC})$ is the algebra of polynomials $\CC[x_1,\dots,x_n]$ 
and for every $p\in\CC[x_1,\dots,x_n]$ we have $p^\perp(x_1,\dots,x_n)=p(-x_1,\dots,-x_n)$. 
\qed
\end{example}

We are going to describe in Remark~\ref{univ_rem} 
some realizations of the universal enveloping algebra $\U(\gg_{\CC})$ 
which are needed later. 
To this end we first introduce the regular representations of $G$ on distribution spaces. 

\begin{definition}\label{reg}
\normalfont
We shall need the \emph{translation maps} $\lambda_g\colon G\to G$, $x\mapsto gx$ 
and $\rho_g\colon G\to G$, $x\mapsto xg$
defined by an arbitrary element $g\in G$. 
The corresponding \emph{regular representations} of $G$ on the distribution space $\Dc'(G)$ 
are defined by 
$$\lambda\colon G\to\End(\Dc'(G)),\quad \langle\lambda(g)u,\phi\rangle=\langle u,\phi\circ\lambda_g\rangle$$
and 
$$\rho\colon G\to\End(\Dc'(G)),\quad \langle\rho(g)u,\phi\rangle=\langle u,\phi\circ\rho_{g^{-1}}\rangle$$
whenever $g\in G$, $u\in\Dc'(G)$, and $\phi\in\Dc(G)$. 
For every $X\in\gg$ and $\phi\in\Ec(G)$ we also define the functions 
$$\dot\lambda(X)\phi\colon G\to\CC, \quad 
(\dot\lambda(X)\phi)(z)=\frac{\de}{\de t}\Big\vert_{t=0}\phi(\exp_G(-tX)z)$$
and 
$$\dot\rho(X)\phi\colon G\to\CC, \quad 
(\dot\rho(X)\phi)(z)=\frac{\de}{\de t}\Big\vert_{t=0}\phi(z\exp_G(tX)).$$
Then we can define the \emph{derivatives of the regular representations} by 
$$\dot\lambda\colon\gg\to\End(\Dc'(G)),\quad 
\langle\dot\lambda(X)u,\phi\rangle:=\langle u,\dot\lambda(-X)\phi\rangle$$
and 
$$\dot\rho\colon\gg\to\End(\Dc'(G)),\quad 
\langle\dot\rho(X)u,\phi\rangle:=\langle u,\dot\rho(-X)\phi\rangle.$$
These derivatives are homomorphisms of Lie algebras, 
hence condition~\eqref{univ_def_item2} in Definition~\ref{univ_def}  
shows that they can be uniquely extended to unital homomorphisms 
of associative algebras $\U(\gg_{\CC})\to\End(\Dc'(G))$. 
These extensions will also be denoted by 
$\dot\lambda\colon\U(\gg_{\CC})\to\End(\Dc'(G))$ and 
$\dot\rho\colon\U(\gg_{\CC})\to\End(\Dc'(G))$, respectively. 

For later use, we also introduce the notation $\phi^\perp(x):=\phi(x^{-1})$ 
for every $\phi\in\Ec(G)$ and $x\in G$. 
This gives rise to the idempotent linear mapping 
$$\Dc'(G)\to\Dc'(G),\quad u\mapsto u^\perp,$$
where $\langle u^\perp,\phi\rangle:=\langle u,\phi^\perp\rangle$ 
for $u\in\Dc'(G)$ and $\phi\in\Dc(G)$. 
\qed
\end{definition}

\begin{remark}\label{univ_rem} 
\normalfont
With Definition~\ref{reg} at hand, 
we can describe some realizations of the universal enveloping algebra $\U(\gg_{\CC})$ as follows. 
For the sake of simplicity, let us denote by $\Ec'_{\1}(G)$ the space of distributions on $G$ 
with the support contained in $\{\1\}$, thought of as a complex unital associative algebra 
with respect to the convolution product, cf.~Example~\ref{conv_ex}.  
(This set should actually be denoted by $\Ec'_{\{\1\}}(G)$ according to Remark~\ref{comp_supp}.) 
Recall that $\delta_{\1}\in\Ec'(G)$ is the Dirac distribution at $\1\in G$. 

Both mappings
$$\begin{aligned} 
\U(\gg_{\CC})\to\Ec'_{\1}(G), \quad & w\mapsto\dot\lambda(w)\delta_{\1}, \\
\U(\gg_{\CC})\to\Ec'_{\1}(G),\quad & w\mapsto\dot\rho(w)\delta_{\1}
  \end{aligned}$$
are isomorphisms of complex unital associative algebras.  
(See for instance Th.~1 in Sect.~10.4 of~\cite{Ki76}.)
These isomorphisms are related by the commutative diagram 
$$\xymatrix{\U(\gg_{\CC}) \ar[r]^{w\mapsto w^\perp} \ar[d]_{\dot\lambda(\cdot)\delta_{\1}} 
      & \U(\gg_{\CC})\ar[d]^{\dot\rho(\cdot)\delta_{\1}} \\
\Ec'_{\1}(G) \ar[r]^{u\mapsto u^\perp} & \Ec'_{\1}(G) 
}$$
where the horizontal arrows stand for 
the mappings introduced in Definitions \ref{univ_def} and~\ref{reg}, 
respectively. 
From now on, \emph{we perform the identification $\U(\gg_{\CC})\simeq\Ec'_{\1}(G)$ by means of the mapping 
$w\mapsto\dot\rho(w)\delta_{\1}$}, by writing simply $w$ instead of $\dot\rho(w)\delta_{\1}$ whenever $w\in\U(\gg_{\CC})$. 
\qed
\end{remark}

\begin{proposition}\label{dirac_dec}
For every integer $m\ge1$ and every compact neighbourhood $K$ of $\1\in G$ 
there exist finitely many elements $u_1,\dots,u_N\in\U(\gg_{\CC})$ 
and the functions $\phi_1,\dots,\phi_N\in\Cc^m_0(G)$ with 
$\bigcup\limits_{j=1}^N\supp\phi_j\subseteq K$ such that 
$\delta_{\1}=\sum\limits_{j=1}^N \phi_j\ast u_j$. 
\end{proposition}

\begin{proof}
Use Lemme~2 in \cite{Ro76} or Lemma~2.3 in \cite{DDJP09}; 
see also the proof of Lemme~1.1 in~\cite{Ca76}. 
\end{proof}

\subsection{Smooth vectors} 

\begin{definition}\label{deriv_def}
\normalfont
The  \emph{smooth vectors} for 
the representation $\pi\colon G\to\Bc(\Yc)$ 
are the elements of the linear subspace of $\Yc$ defined by 
$$\Yc_\infty:=\{y\in\Yc\mid\pi(\cdot)y\in\Ec(G,\Yc)\}. $$
The linear space $\Yc_\infty$ will be endowed with the linear topology 
which makes the linear injective map 
$$\Yc_\infty\to\Ec(G,\Yc),\quad y\mapsto\pi(\cdot)y $$
into a linear topological isomorphism onto its image. 

For every distribution $u\in\Ec'(G)$ and 
every smooth vector $y\in\Yc_\infty$ we then define 
$$\dot\pi(u)y:=\langle u\hotimes\id_{\Yc},\pi(\cdot)y\rangle\in\Yc.$$
\qed
\end{definition}

\begin{proposition}\label{deriv_prop}
The following assertions hold: 
\begin{enumerate}
\item\label{deriv_prop_item0} 
The space of smooth vectors $\Yc_\infty$ is a Fr\'echet space 
and the inclusion map $\Yc_\infty\hookrightarrow\Yc$ is continuous. 
\item\label{deriv_prop_item1} 
The space $\Yc_\infty$ is dense in $\Yc$. 
\item\label{deriv_prop_item2} 
For every $u\in\Ec'(G)$ we have $\dot\pi(u)\Yc_\infty\subseteq\Yc_\infty$. 
\item\label{deriv_prop_item3} 
The mapping $\dot\pi\colon\Ec'(G)\to\End(\Yc_\infty)$ 
is a homomorphism of unital associative algebras. 
\item\label{deriv_prop_item4}  
For every $X\in\gg$ and $y\in\Yc_\infty$ 
we have $\dot\pi(X)y:=\frac{\de}{\de t}\Big\vert_{t=0}\pi(\exp_G(tX))y$.
\item\label{deriv_prop_item5} 
For every $y\in\Yc_\infty$ 
we have a smooth mapping $\pi(\cdot)y\colon G\to\Yc_\infty$. 
\end{enumerate}
\end{proposition}

\begin{proof}
See for instance \cite{Wa72} and Sect.~10.5 in \cite{Ki76}.
\end{proof}

\begin{notation}
\normalfont 
We shall always denote by $\Yc_{-\infty}$ the space of continuous antilinear functionals on 
the Fr\'echet space $\Yc_\infty$. 
\qed
\end{notation}

\begin{proposition}\label{dm}
For every integer $m\ge1$ there exist finitely many functions $\phi_1,\dots,\phi_N\in\Cc^m_0(G)$ 
such that for every $y\in\Yc_\infty$ there exist $y_1,\dots,y_N\in\Yc$ 
satisfying the equality $y=\dot\pi(\phi_1)y_1+\cdots+\dot\pi(\phi_N)y_N$. 
\end{proposition}

\begin{proof}
Use Proposition~\ref{dirac_dec} to get  $u_1,\dots,u_N\in\U(\gg_{\CC})$ 
and $\phi_1,\dots,\phi_N\in\Cc^m_0(G)$ with 
$\delta_{\1}=\sum\limits_{j=1}^N \phi_j\ast u_j$. 
Then Proposition~\ref{deriv_prop} shows that 
$$y=\dot\pi(\delta_{\1})y=\sum_{j=1}^N \dot\pi(\phi_j)\dot\pi(u_j)y
=\sum_{j=1}^N \dot\pi(\phi_j)y_j,$$
where we have denoted $y_j:=\dot\pi(u_j)y$ for $j=1,\dots,N$.
\end{proof}

\begin{remark}
\normalfont
As we already mentioned, the smooth vectors for representations of Lie groups 
were discussed in detail in \cite{Wa72}. 
Other important references in this connection are \cite{Ga47}, \cite{Ga60}, \cite{Ki76}, \cite{Ca76}, \cite{CGP77}, \cite{DM78}, \cite{RT87}, and \cite{CG90}.
\qed
\end{remark}

\section{Smooth operators for unitary representations}

We are going to discuss here the space of smooth operators for 
a given representation of a Lie group. 
The method of investigation was suggested in \cite{Pe94} 
and relies on exhibiting 
this space of operators as the space of smooth vectors 
for a suitable representation (see Definition~\ref{suitable}). 
The main result is recorded as Theorem~\ref{smooth_ideal} 
and it is particularly significant in the case of unitary irreducible representations of nilpotent Lie groups (Corollary~\ref{smooth_ideal_cor}).

\begin{notation}
\normalfont
In this section we shall use the following notation:
\begin{itemize}
\item $G$ is a connected unimodular Lie group with the Lie algebra $\gg$; 
\item $\de x$ denotes a fixed Haar measure on $G$; 
\item $\Hc$ is a complex Hilbert space; 
\item $\pi\colon G\to \Bc(\Hc)$ is a continuous unitary representation, 
and $\Hc_\infty$ is the corresponding space of smooth vectors.
\end{itemize}
\qed
\end{notation}

The following notion of smooth operator 
was singled out on page~349 in \cite{How77} 
and then further developed in \cite{Pe94}. 

\begin{definition}\label{smooth_oper}
\normalfont
The set $\Bc(\Hc)_\infty$ of \emph{smooth operators} for the representation~$\pi$ 
is defined as the set of all operators $T\in\Bc(\Hc)$ satisfying the following conditions: 
\begin{enumerate}
\item\label{smooth_oper_item1} 
We have $T(\Hc)+T^*(\Hc)\subseteq\Hc_\infty$. 
\item\label{smooth_oper_item2} 
For every $u\in\U(\gg_{\CC})$ the operators 
$\dot\pi(u)T$ and $\dot\pi(u)T^*$ are bounded on~$\Hc$. 
\end{enumerate}
\qed
\end{definition}

\begin{example}\label{smooth_oper_ex}
\normalfont
For every $x,y\in\Hc_\infty$ the rank-one operator $(\cdot\mid x)y$ belongs to 
the space of smooth operators~$\Bc(\Hc)_\infty$. 
We shall see in Corollary~\ref{approximation} that 
the linear span of these rank-one operators is dense in $\Bc(\Hc)_\infty$ 
provided that $G$ is a nilpotent Lie group and $\pi$ is an irreducible representation. 
\qed
\end{example}

\begin{remark}\label{smooth_oper_rem}
\normalfont
It follows at once by Definition~\ref{smooth_oper} that $\Bc(\Hc)_\infty$ is an associative $*$-subalgebra 
of $\Bc(\Hc)$. 
\qed
\end{remark}

\begin{definition}\label{char_def}
\normalfont
We shall say that the representation $\pi$ \emph{has a smooth character} if 
for every $\phi\in\Dc(G)$ we have $\pi(\phi)\in\Sg_1(\Hc)$ 
and the linear mapping 
$$\Dc(G)\to\Sg_1(\Hc),\quad \phi\mapsto\dot\pi(\phi) $$
is continuous. 
In this case we define the corresponding \emph{character} as
$$\chi_\pi\colon\Dc(G)\to\CC,\quad \chi_\pi(\phi):=\Tr\pi(\phi).$$
Note that $\chi_\pi\in\Dc'(G)$. 
\qed
\end{definition}

\begin{example}
\normalfont
Every unitary irreducible representation of a nilpotent Lie group has a smooth character; 
see for instance Th.~2 in \S5 of Ch.~II, Part.II in \cite{Pu67}.
\qed
\end{example}

\begin{remark}\label{char_rem}
\normalfont
If the representation $\pi$ has a smooth character, then there exists a continuous seminorm $p(\cdot)$ 
on $\Dc(G)$ such that 
$$(\forall \phi\in\Dc(G))\quad \Vert\dot\pi(\phi)\Vert_1\le p(\phi).$$
In view of the definition of the topology on $\Dc(G)$ and of 
the fact that $\Dc(G)$ is dense in $\Cc^m_0(G)$ for every $m\ge1$, 
it then easily follows that for every compact subset $K\subset G$ 
there exists an integer $m\ge1$ 
such that for every $\phi\in\Cc^m_0(G)\cap\Ec'_K(G)$ we have $\dot\pi(\phi)\in\Sg_1(\Hc)$, 
and moreover the mapping 
$$\Cc^m_0(G)\cap\Ec'_K(G)\to\Sg_1(\Hc),\quad \phi\mapsto\dot\pi(\phi) $$
is linear and continuous. 
\qed
\end{remark}

\begin{definition}\label{adm}
\normalfont
An \emph{admissible ideal} is a non-trivial two-sided ideal $\Jc$ of $\Bc(\Hc)$ 
with the following properties: 
\begin{enumerate}
\item\label{adm_item1} 
The ideal $\Jc$ is endowed with a complete norm $\Vert\cdot\Vert_{\Jc}$ such that 
for every $A,B\in\Bc(\Hc)$ and $T\in\Jc$ we have 
$\Vert ATB\Vert_{\Jc}\le\Vert A\Vert \cdot\Vert T\Vert_{\Jc}\cdot\Vert B\Vert$ 
and $\Vert T^*\Vert_{\Jc}=\Vert T\Vert_{\Jc}$. 
\item\label{adm_item2} 
The ideal $\Fc(\Hc)$ of finite-rank operators is a dense subspace of $\Jc$. 
\item\label{adm_item3} 
For every $x,y\in\Hc$ we have 
$\Vert(\cdot\mid x)y\Vert_{\Jc}=\Vert x\Vert\cdot\Vert y\Vert$.
\end{enumerate}
\qed
\end{definition}

\begin{example}
\normalfont 
Every Schatten ideal $\Sg_p(\Hc)$ with $1\le p\le\infty$ is an admissible ideal. 
There exist many other examples of admissible ideals; see for instance \cite{GK69}. 
\qed 
\end{example}

\begin{remark}\label{so}
\normalfont 
Let $\Jc$ be an admissible ideal. 
By using condition~\eqref{adm_item2} in Definition~\ref{adm} 
with $A=\id_{\Hc}$ and $B=(\cdot\mid x)x$ for $x\in\Hc$, 
and then taking into account 
condition~\eqref{adm_item3}, it follows that 
$\Vert Tx\Vert\le\Vert T\Vert_{\Jc}\cdot\Vert x\Vert$. 
That is, for every $T\in\Jc$ we have $\Vert T\Vert\le\Vert T\Vert_{\Jc}$. 

On the other hand, it follows at once by condition~\eqref{adm_item3} in Definition~\ref{adm} that for every $T\in\Fc(\Hc)$ we have 
$\Vert T\Vert_{\Jc}\le\Vert T\Vert_1$. 
Since $\Fc(\Hc)$ is dense in $\Sg_1(\Hc)$, 
we get 
$$(\forall T\in\Sg_1(\Hc))\quad 
\Vert T\Vert\le\Vert T\Vert_{\Jc}\le\Vert T\Vert_1. $$
In particular, we have $\Sg_1(\Hc)\subseteq\Jc$. 
\qed
\end{remark}

\begin{definition}\label{suitable}
\normalfont 
For every admissible ideal $\Jc$ we define 
a linear representation $\pi^{\otimes 2}_{\Jc}\colon G\times G\to\Bc(\Jc)$ by 
$$\pi^{\otimes 2}_{\Jc}(g_1,g_2)T:=\pi(g_1)T\pi(g_2)^{-1} $$
for every $g_1,g_2\in G$ and $T\in\Jc$. 
\qed
\end{definition}

\begin{lemma}
The representation $\pi^{\otimes 2}_{\Jc}\colon G\times G\to\Bc(\Jc)$ 
is continuous for every admissible ideal $\Jc$. 
\end{lemma}

\begin{proof}
If $T\in\Fc(\Hc)$, then it is straightforward to check that 
$\pi^{\otimes 2}_{\Jc}(\cdot)T\colon G\times G\to\Jc$ is a continuous mapping. 

Now let $T\in\Jc$ arbitrary. 
Since $\Jc$ is admissible, there exists a sequence $\{T_k\}_{k\ge1}$ in $\Fc(\Hc)$ 
such that $\lim\limits_{k\to\infty}\Vert T-T_k\Vert_{\Jc}=0$. 
On the other hand, since $\pi$ is a unitary representation, it follows that 
for $k=1,2,\dots$ and every $(g_1,g_2)\in G\times G$ we have 
$$\Vert\pi^{\otimes 2}_{\Jc}(g_1,g_2)T-\pi^{\otimes 2}_{\Jc}(g_1,g_2)T_k\Vert_{\Jc}\le\Vert T-T_k\Vert_{\Jc}.$$
Therefore $\pi^{\otimes 2}_{\Jc}(\cdot)T\colon G\times G\to\Jc$ is the uniform limit on $G\times G$ 
of the sequence of continuous mappings $\pi^{\otimes 2}_{\Jc}(\cdot)T$, 
hence it is in turn continuous. 
\end{proof}

\begin{theorem}\label{smooth_ideal}
Let $\pi\colon G\to\Bc(\Hc)$ be a continuous unitary representation, 
assume that $\Jc\subset\Bc(\Hc)$ is an admissible ideal, 
and denote by $\Jc_\infty$ the space of smooth vectors 
for the corresponding representation~$\pi^{\otimes 2}_{\Jc}$. 
Then the following assertions hold: 
\begin{enumerate}
\item\label{smooth_ideal_item1} 
We have $\Jc_\infty\subseteq\Bc(\Hc)_\infty$.
\item\label{smooth_ideal_item2} 
If the Fr\'echet space of smooth vectors $\Hc_\infty$ is nuclear, 
then we have 
$$\Jc_\infty=\Bc(\Hc)_\infty\subseteq\Sg_1(\Hc),$$  
and the Fr\'echet space $\Jc_\infty$ does not depend on 
the choice of the admissible ideal~$\Jc$.
\end{enumerate}
\end{theorem}

\begin{proof}
\eqref{smooth_ideal_item1}
To prove the inclusion $\Jc_\infty\subseteq\Bc(\Hc)_\infty$,  
let $T\in\Jc_\infty$ arbitrary,    
hence the mapping 
$$G\times G\to\Jc,\quad 
(g_1,g_2)\mapsto\pi^{\otimes 2}_{\Jc}(g_1,g_2)T=\pi(g_1)T\pi(g_2)^{-1}$$
is smooth.
In particular, 
the mapping $\pi(\cdot)T\colon G\to\Jc$ is smooth. 
On the other hand, it follows by Remark~\ref{so} that 
for arbitrary $x\in\Hc$ we have a continuous linear mapping 
$\Jc\to\Hc$, $T\mapsto Tx$. 
Hence the mapping $\pi(\cdot)Tx\colon G\to\Hc$ will be smooth as 
a composition of two smooth mappings. 
Thus for arbitrary $x\in\Hc$ we have $Tx\in\Hc_\infty$. 
Moreover, since the operation of taking the Hilbert space adjoint 
is ($\RR$-linear and) continuous on $\Jc$ by condition 
\eqref{adm_item1} in Definition~\ref{adm}, 
it follows at once that $T^*\in\Jc_\infty$. 
Hence by the above reasoning with $T$ replaced by $T^*$ 
we get $T^*x\in\Hc_\infty$ for arbitrary $x\in\Hc$. 
Thus the operator $T$ satisfies condition~\eqref{smooth_oper_item1} 
in Definition~\ref{smooth_oper}. 
To check condition~\eqref{smooth_oper_item2} in the same definition 
just note that since the mapping $\pi(\cdot)T\colon G\to\Jc$ is smooth, 
it follows that for every $u\in\U(\gg_{\CC})$ we have 
$\dot{\pi}(u)T\in\Jc$, hence $\dot{\pi}(u)T\in\Bc(\Hc)$. 
Since we have seen above that $T^*\in\Jc_\infty$, 
it also follows that $\dot{\pi}(u)T^*\in\Bc(\Hc)$. 
This completes the proof of the fact that $T\in\Bc(\Hc)_\infty$. 

\eqref{smooth_ideal_item2} 
If $\Hc_\infty$ is a nuclear space, 
then the inclusion map $\Hc_\infty\hookrightarrow\Hc$ is a nuclear operator
(see Prop.~7.2 in Ch.~III of \cite{Sch66}). 
Since condition~\eqref{smooth_oper_item1} in Definition~\ref{smooth_oper} 
shows an arbitrary operator $T\in\Bc(\Hc)_\infty$ 
factorizes as 
$\Hc\mathop{\longrightarrow}\limits^T\Hc_\infty\hookrightarrow\Hc$, 
it follows that $T\in\Sg_1(\Hc)$ 
(see Cor.~2 to Prop.~7.2 in Ch.~III of \cite{Sch66}). 
Thus, by taking into account the above Assertion~\eqref{smooth_ideal_item1}, 
we get  
\begin{equation}\label{smooth_ideal_eq1}
\Jc_\infty\subseteq\Bc(\Hc)_\infty\subseteq\Sg_1(\Hc).
\end{equation}
To see that $\Jc_\infty$ does not depend on the choice of 
the admissible ideal $\Jc$, we shall prove the equality of Fr\'echet spaces
\begin{equation}\label{smooth_ideal_eq2}
\Jc_\infty=\Sg_1(\Hc)_\infty, 
\end{equation}
where the right-hand side denotes the space of smooth vectors for 
the representation 
$\pi^{\otimes 2}_{\Sg_1(\Hc)}\colon G\times G\to\Bc(\Sg_1(\Hc))$. 
First recall from Remark~\ref{so} that 
we have a continuous inclusion map $\Sg_1(\Hc)\hookrightarrow\Jc$, 
which clearly intertwines the representations 
$\pi^{\otimes 2}_{\Sg_1(\Hc)}$ and $\pi^{\otimes 2}_{\Jc}$. 
It then easily follows by Definition~\ref{deriv_def} that 
we have a continuous inclusion map for the corresponding spaces of smooth vectors $\Sg_1(\Hc)_\infty\hookrightarrow\Jc_\infty$. 
On the other hand, we have already proved that $\Jc_\infty\subseteq\Sg_1(\Hc)$, 
hence there exists the following  commutative diagram 
$$\xymatrix{
 \Jc_\infty\ar[r] \ar[d] & \Sg_1(\Hc)\ar[ld]   \\
  \Jc              }
$$
whose arrows stand for inclusion maps between Fr\'echet or Banach spaces. 
The arrows that point to $\Jc$ are continuous inclusions 
(by Remark~\ref{so} and Proposition~\ref{deriv_prop}\eqref{deriv_prop_item0}), 
hence the closed graph theorem implies that 
the inclusion map $\Jc_\infty\hookrightarrow\Sg_1(\Hc)$ is continuous as well. 
Since for arbitrary $T\in\Jc_\infty$ 
the mapping $G\times G\to\Jc_\infty$, $(g_1,g_2)\mapsto\pi(g_1)T\pi(g_2)^{-1}$ 
is smooth by Proposition~\ref{deriv_prop}\eqref{deriv_prop_item5}, 
it then follows that the mapping 
$G\times G\to\Sg_1(\Hc)$, $(g_1,g_2)\mapsto\pi(g_1)T\pi(g_2)^{-1}$  
is also smooth, hence $T\in\Sg_1(\Hc)_\infty$. 
Thus $\Jc_\infty=\Sg_1(\Hc)_\infty$ as sets. 
Since both sides of this equality are Fr\'echet spaces and we have already seen 
that the inclusion map $\Sg_1(\Hc)_\infty\hookrightarrow\Jc_\infty$ 
is continuous, it follows by the open mapping theorem that we have 
the equality of Fr\'echet spaces in \eqref{smooth_ideal_eq2}. 

Finally, note that for arbitrary $T\in\Bc(\Hc)_\infty$ and 
every $u\in\U(\gg_{\CC})$ we have $\dot\pi(u)T,\dot\pi(u)T^*\in\Bc(\Hc)_\infty$ 
(see Definition~\ref{smooth_oper}). 
On the other hand, 
we have proved above that $\Bc(\Hc)_\infty\subseteq\Sg_1(\Hc)$, 
hence $\dot\pi(u)T,\dot\pi(u)T^*\in\Sg_1(\Hc)$ for all $u\in\U(\gg_{\CC})$, 
and this implies that $T\in\Sg_1(\Hc)_\infty$. 
Thus $\Bc(\Hc)_\infty\subseteq\Sg_1(\Hc)_\infty$, 
and then by using Assertion~\eqref{smooth_ideal_item1} 
with $\Jc=\Sg_1(\Hc)$ we get $\Bc(\Hc)_\infty=\Sg_1(\Hc)_\infty$. 
Now by \eqref{smooth_ideal_eq2} and \eqref{smooth_ideal_eq1} 
we get $\Jc_\infty=\Bc(\Hc)_\infty\subseteq\Sg_1(\Hc)$, 
and this completes the proof. 
\end{proof}

\begin{corollary}\label{smooth_ideal_cor}
Assume that $G$ is a nilpotent Lie group and 
$\pi\colon G\to\Bc(\Hc)$ is a unitary irreducible representation.  
If $\Jc\subset\Bc(\Hc)$ is an admissible ideal, 
and we denote by $\Jc_\infty$ the space of smooth vectors 
for the corresponding representation~$\pi^{\otimes 2}_{\Jc}$,  
Then the following assertions hold: 
\begin{enumerate}
\item\label{smooth_ideal_cor_item1} 
The Fr\'echet space of smooth vectors $\Hc_\infty$ is nuclear. 
\item\label{smooth_ideal_cor_item2} 
We have 
$\Jc_\infty=\Bc(\Hc)_\infty\subseteq\Sg_1(\Hc)$,   
and the Fr\'echet space $\Jc_\infty$ does not depend on 
the choice of the admissible ideal~$\Jc$.
\item\label{smooth_ideal_cor_item3} 
The space of smooth operators $\Bc(\Hc)_\infty$ has the natural structure of  
a nuclear Fr\'echet space.
\end{enumerate}
\end{corollary}

\begin{proof}
Since the representation $\pi$ is irreducible, 
there exists a linear topological isomorphism from 
the Fr\'echet space $\Hc_\infty$ onto 
the Schwartz space of rapidly decreasing functions $\Sc(\RR^{d/2})$, 
where $d$ is equal to the dimension of the coadjoint orbit 
of $G$ corresponding to the representation~$\pi$.  
(This follows by Th.~1 in \S 5 of Ch.~II, Part.~II in \cite{Pu67}; 
see also the Cor.~to Th.~3.1 in \cite{CGP77}, or \cite{CG90}.) 
On the other hand, it is well known that 
the Schwartz space $\Sc(\RR^{d/2})$ 
is nuclear (see for instance Ex.~5 in \S 8 of Ch.~III in \cite{Sch66}). 
Therefore the Fr\'echet space $\Hc_\infty$ is nuclear, 
and then Theorem~\ref{smooth_ideal} applies. 

Finally, by using Assertion~\eqref{smooth_ideal_cor_item2} 
when $\Jc=\Sg_2(\Hc)$ (the Hilbert-Schmidt ideal), it follows that 
$\Bc(\Hc)_\infty$ is equal to the space of smooth vectors 
for the unitary representation $\pi^{\otimes 2}_{\Sg_2(\Hc)}$, 
hence it is a Fr\'echet space in a natural way. 
On the other hand, the representation 
$\pi^{\otimes 2}_{\Sg_2(\Hc)}$ 
is irreducible since so is~$\pi$. 
(See for instance the proof of Lemma~2.18(a) in \cite{BB09c}.) 
Now the fact that $\Bc(\Hc)_\infty$ is nuclear follows 
by the above Assertion~\eqref{smooth_ideal_cor_item1} applied 
for the unitary irreducible 
representation~$\pi^{\otimes 2}_{\Sg_2(\Hc)}\colon G\times G\to\Bc(\Sg_2(\Hc))$. 
\end{proof}

\begin{corollary}\label{regularization}
If $G$ is a nilpotent Lie group and 
$\pi\colon G\to\Bc(\Hc)$ is a unitary irreducible representation, 
then the operators in $\Bc(\Hc)_\infty$ are precisely the regularizing operators. 
That is, $A\in\Bc(\Hc)_\infty$  
if and only if $A$ extends to a continuous linear map $A\colon\Hc_{-\infty}\to\Hc_\infty$, 
so that the diagram   
$$\xymatrix{
\Hc_{-\infty} \ar[r]^{A} & \Hc_\infty \ar@{^{(}->}[d] \\
\Hc \ar[r]^A \ar@{^{(}->}[u] & \Hc
} $$
is commutative
\end{corollary}

\begin{proof}
By the closed graph theorem, it is enough to prove that 
if $A\in\Bc(\Hc)_\infty$  and $f\in\Hc_{-\infty}$,  
then $Af\in\Hc_\infty$, 
in the sense that there exists a smooth vector 
denoted $Af$ such that for every $\phi\in\Hc_\infty$ we have 
$(f\mid A^*\phi)=(Af\mid\phi)$. 
This is a consequence of the above Proposition~\ref{deriv_prop}\eqref{deriv_prop_item2}, Corollary\ref{smooth_ideal_cor}, and  Th.~1.3(b) in~\cite{Ca76}. 

Conversely, it follows by Definition~\ref{smooth_oper} 
that the restriction to $\Hc$ of every continuous linear map 
$A\colon\Hc_{-\infty}\to\Hc_\infty$ 
belongs to $\Bc(\Hc)_\infty$.
\end{proof}

\begin{corollary}\label{approximation}
Assume that $G$ is a nilpotent Lie group and 
$\pi\colon G\to\Bc(\Hc)$ is a unitary irreducible representation.  
The linear space spanned by the operators $(\cdot\mid x)y$ with $x,y\in\Hc_\infty$ 
is dense in $\Bc(\Hc)_\infty$. 
\end{corollary}

\begin{proof}
Let $T\in\Bc(\Hc)_\infty$ arbitrary. 
Then Corollary~\ref{smooth_ideal_cor}\eqref{smooth_ideal_cor_item2} 
shows that $T$ is a smooth vector for the representation 
$\pi^{\otimes 2}_{\Sg_2(\Hc)}\colon G\times G\to\Bc(\Sg_2(\Hc))$. 
It follows by Proposition~\ref{dm} that 
there exist finitely many functions 
$\phi_1,\dots,\phi_N\in\Cc^m_0(G\times G)$ 
and the operators $Y_1,\dots,Y_N\in\Sg_2(\Hc)$ such that 
$$T=\pi^{\otimes 2}_{\Sg_2(\Hc)}(\phi_1)Y_1+\cdots+\pi^{\otimes 2}_{\Sg_2(\Hc)}(\phi_N)Y_N.$$
Since $\Dc(G\times G)$ is dense in $\Cc^m_0(G\times G)$ 
and $\Dc(G)\otimes\Dc(G)$ is dense in $\Dc(G\times G)$, 
it follows by Proposition~\ref{deriv_prop}\eqref{deriv_prop_item3} 
that $T$ can be approximated in $\Bc(\Hc)_\infty$ by finite linear combinations of operators of the form 
$$\pi^{\otimes 2}_{\Sg_2(\Hc)}(\psi_1\otimes\psi_2)Y=\pi(\psi_1)Y\pi(\psi_2^\perp)$$
with $\psi_1,\psi_2\in\Dc(G)$ and $Y\in\Sg_2(\Hc)$. 
On the other hand, such an $Y$ can be approximated in $\Sg_2(\Hc)$ 
by finite linear combinations of operators $(\cdot\mid v_2)v_1$ 
with $v_1,v_2\in\Hc$. 
The corollary now follows by noticing that 
$$\pi(\psi_1)((\cdot\mid v_2)v_1)\pi(\psi_2^\perp)=(\cdot\mid \pi(\bar\psi_2^\perp)v_2)\pi(\psi_1)v_1$$
and recalling that $\pi(\psi)v\in\Hc_\infty$ when $\psi\in\Dc(G)$ and $v\in\Hc$ 
(\cite{Ga47}).
\end{proof}

\begin{remark}\label{smooth_ideal_rem}
\normalfont 
Let $\Jc\subset\Bc(\Hc)$ be any admissible ideal. 
If the representation $\pi$ has a smooth character, 
then the corresponding space of smooth vectors $\Hc_\infty$ is nuclear 
according to Th.~2.6 in \cite{Ca76}, 
hence the above Theorem~\ref{smooth_ideal}\eqref{smooth_ideal_item2} applies. 

The inclusion $\Jc_\infty\subseteq\Sg_1(\Hc)$ 
can be alternatively proved in this case as follows.  
Let $T\in\Jc_\infty$ arbitrary. 
Since the representation $\pi$ has a smooth character, 
we have a continuous linear mapping 
$$\Dc(G)\to\Sg_1(\Hc),\quad \phi\mapsto\dot\pi(\phi).$$
On the other hand, note that for every $\phi_1,\phi_2\in\Dc(G)$ we have 
$\pi^{\otimes 2}_{\Jc}(\phi_1\otimes\phi_2)T=\dot\pi(\phi_1)T\dot\pi(\phi_2^\perp)$.  
Hence for arbitrary $Y\in\Jc$ we get a jointly continuous trilinear mapping 
$$\Dc(G)\times\Dc(G)\times\Jc\to\Sg_1(\Hc),\quad 
(\phi_1,\phi_2,Y)\mapsto\pi^{\otimes 2}_{\Jc}(\phi_1\otimes\phi_2)Y, $$
which extends to a jointly continuous bilinear mapping 
$$\Dc(G\times G)\times\Jc\to\Sg_1(\Hc), \quad 
(\phi,Y)\mapsto\pi^{\otimes 2}_{\Jc}(\phi)Y.$$
(See also Remark~\ref{tens}.)
By an argument similar to the one of Remark~\ref{char_rem}, 
we can further extend the above mapping to a continuous bilinear mapping 
\begin{equation}\label{smooth_ideal_rem_eq1}
(\Cc^m_0(G\times G)\cap\Ec'_K(G\times G))\times\Jc\to\Sg_1(\Hc), \quad 
(\phi,Y)\mapsto\pi^{\otimes 2}_{\Jc}(\phi)Y
\end{equation}
for a suitable compact neighborhood $K$ of 
$(\1,\1)\in G\times G$ and a suitably large integer $m\ge1$. 
On the other hand, since $T\in\Jc_\infty$, 
it follows by Proposition~\ref{dm} that 
there exist finitely many functions 
$\phi_1,\dots,\phi_N\in\Cc^m_0(G\times G)\cap\Ec'_K(G\times G)$ 
and the operators $Y_1,\dots,Y_N\in\Jc$ such that 
$T=\pi^{\otimes 2}_{\Jc}(\phi_1)Y_1+\cdots+\pi^{\otimes 2}_{\Jc}(\phi_N)Y_N$, 
hence by \eqref{smooth_ideal_rem_eq1} we get $T\in\Sg_1(\Hc)$. 
\qed
\end{remark}

\section{Weyl-Pedersen calculus}

In the present section we provide a brief discussion 
of the remarkable Weyl correspondence constructed in \cite{Pe94} 
and we shall also describe some complementary results 
which were recently obtained in \cite{BB09c}. 

\subsection{Preduals for coadjoint orbits}
This subsection records some properties of the coadjoint orbits 
of nilpotent Lie groups which play a crucial role for 
the construction of the Weyl-Pedersen calculus. 

\begin{setting}\label{predual_sett}
\normalfont 
We shall use the following notation:
\begin{enumerate}
 \item Let $G$ be a connected, simply connected, nilpotent Lie group with Lie algebra~$\gg$.  
 Then the exponential map $\exp_G\colon\gg\to G$ is a diffeomorphism 
 with the inverse denoted by $\log_G\colon G\to\gg$. 
 \item We denote by $\gg^*$ the linear dual space to $\gg$ and 
  by $\hake{\cdot,\cdot}\colon\gg^*\times\gg\to{\RR}$ the natural duality pairing. 
 \item Let $\xi_0\in\gg^*$ with the corresponding coadjoint orbit $\Oc:=\Ad_G^*(G)\xi_0\subseteq\gg^*$.  
 \item Let $\pi\colon G\to\Bc(\Hc)$ be any unitary irreducible representations 
associated with the coadjoint orbit $\Oc$ by Kirillov's theorem (\cite{Ki62}).
 \item The \emph{isotropy group} at $\xi_0$ is $G_{\xi_0}:=\{g\in G\mid\Ad_G^*(g)\xi_0=\xi_0\}$ 
 with the corresponding \emph{isotropy Lie algebra} $\gg_{\xi_0}=\{X\in\gg\mid\xi_0\circ\ad_{\gg}X=0\}$. 
 If we denote the \emph{center} of $\gg$ by $\zg:=\{X\in\gg\mid[X,\gg]=\{0\}\}$, 
 then it is clear that $\zg\subseteq\gg_{\xi_0}$. 
 \item Let $n:=\dim\gg$ and fix a sequence of ideals in $\gg$, 
$$\{0\}=\gg_0\subset\gg_1\subset\cdots\subset\gg_n=\gg$$
such that $\dim(\gg_j/\gg_{j-1})=1$ and $[\gg,\gg_j]\subseteq\gg_{j-1}$ 
for $j=1,\dots,n$. 
 \item Pick any $X_j\in\gg_j\setminus\gg_{j-1}$ for $j=1,\dots,n$, 
so that the set $\{X_1,\dots,X_n\}$ will be a \emph{Jordan-H\"older basis} in~$\gg$. 
\end{enumerate}
\qed 
\end{setting}

\begin{definition}\label{jump_def}
\normalfont
 Consider the set of \emph{jump indices} of the coadjoint orbit $\Oc$ 
with respect to the aforementioned Jordan-H\"older basis $\{X_1,\dots,X_n\}\subset\gg$, 
$$e:=\{j\in\{1,\dots,n\}\mid \gg_j\not\subseteq\gg_{j-1}+\gg_{\xi_0}\}
=\{j\in\{1,\dots,n\}\mid X_j\not\in\gg_{j-1}+\gg_{\xi_0}\}$$ 
and then define the corresponding \emph{predual of the coadjoint orbit}~$\Oc$, 
$$\gg_e:=\spa\{X_j\mid j\in e\}\subseteq\gg.$$
We note the direct sum decomposition $\gg=\gg_{\xi_0}\dotplus\gg_e$. 
\qed 
\end{definition}

\begin{remark}\label{jump_remark}
\normalfont
Let $\{\xi_1,\dots,\xi_n\}\subset\gg^*$ be the dual basis for $\{X_1,\dots,X_n\}\subset\gg$. 
Then the coadjoint orbit $\Oc$ can be described in terms of the jump indices 
mentioned in Definition~\ref{jump_def}. 
More specifically, if we denote 
$$\gg^*_{\Oc}:=\spa\{\xi_j\mid j\in e\} 
\quad\text{ and }\quad\gg^\perp_{\Oc}:=\spa\{\xi_j\mid j\not\in e\},$$ 
then the coadjoint orbit  
$\Oc\subset\gg^*\simeq\gg_e^*\times\gg_e^\perp$ 
is the graph of a certain \emph{polynomial mapping} $\gg_e^*\to\gg_e^\perp$. 
This leads to the following pieces of information on $\Oc$: 
\begin{enumerate}
 \item $\dim\Oc=\dim\gg_e=\card e=:d$;  
 \item if we let $j_1<\cdots<j_d$ such that $e=\{j_1,\dots,j_d\}$, 
then the mapping 
$$\Oc\to{\RR}^d,\quad \xi\to(\hake{\xi,X_{j_1}},\dots,\hake{\xi,X_{j_d}}) $$
is a global chart which takes the Liouville measure of $\Oc$ 
to a Lebesgue measure on~${\RR}^d$. 
\end{enumerate}
We define the Fourier transform $\Sc(\Oc)\to\Sc(\gg_e)$ by 
$$(\forall X\in\gg_e)\quad \widehat{a}(X)=\int\limits_{\Oc}\ee^{-\ie\hake{\xi,X}}a(\xi)\de\xi  $$
for every $a\in\Sc(\Oc)$, where $\de\xi$ stands for a Liouville measure on $\Oc$. 
This Fourier transform is invertible. 
The Lebesgue measure on $\gg_e$ can be normalized 
such that the Fourier transform extends to a unitary operator 
$$L^2(\Oc)\to L^2(\gg_e),\quad a\mapsto\widehat{a},$$
and its inverse is defined by the usual formula. 
We shall always consider the predual $\gg_e$ endowed with 
this normalized measure. 
(See for instance Lemma~1.6.1 in \cite{Pe89} and Lemma~4.1.1 in \cite{Pe94} 
for more details and proofs for the above assertions.) 
\qed
\end{remark}

\begin{remark}
\normalfont
Some basic references for the geometry of coadjoint orbits 
of nilpotent Lie groups include \cite{Pu67}, \cite{Pe84}, \cite{Pe88}, \cite{Pe89}, 
and \cite{CG90}; see also \cite{BB09d}. 
\qed
\end{remark}

\subsection{Weyl-Pedersen calculus and Moyal identities}
We begin this subsection by the general construction of 
a Weyl correspondence due to \cite{Pe94}. 

\begin{definition}\label{calc_def}
\normalfont 
The \emph{Weyl-Pedersen calculus} $\Op^\pi(\cdot)$ for the unitary representation~$\pi$ is defined 
for every $a\in\Sc(\Oc)$ by 
$$\Op^\pi(a)=\int\limits_{\gg_e}\widehat{a}(X)\pi(\exp_GX)\de X\in\Bc(\Hc). $$
We call $\Op^\pi(a)$ is the \emph{pseudo-differential operator} with 
the \emph{symbol} $a\in\Sc(\Oc)$. 
\qed 
\end{definition}

\begin{theorem}\label{pedersen}
The Weyl-Pedersen calculus has the following properties: 
\begin{enumerate}
\item\label{pedersen_item1} For every symbol $a\in\Sc(\Oc)$ we have $\Op^\pi(a)\in\Bc(\Hc)_\infty$ 
and the mapping $\Sc(\Oc)\to\Bc(\Hc)_\infty$, $a\mapsto\Op^\pi(a)$ is 
a linear topological isomorphism. 
\item\label{pedersen_item1bis} 
For every $T\in\Bc(\Hc)_\infty$ we have $T=\Op^\pi(a)$, 
where $a\in\Sc(\Oc)$ satisfies the condition $\widehat{a}(X)=\Tr(\pi(\exp_G X)^{-1}T)$ for every $X\in\gg_e$. 
\item\label{pedersen_item2} For every $a,b\in\Sc(\Oc)$ we have 
\begin{enumerate}
\item\label{pedersen_item2a} $\Op^\pi(\bar{a})=\Op^\pi(a)^*$; 
\item\label{pedersen_item2b} $\Tr(\Op^\pi(a))=\int\limits_{\Oc}a(\xi)\de\xi$; 
\item\label{pedersen_item2c} $\Tr(\Op^\pi(a)\Op^\pi(b))=\int\limits_{\Oc}a(\xi)b(\xi)\de\xi$;
\item\label{pedersen_item2d} $\Tr(\Op^\pi(a)\Op^\pi(b)^*)=\int\limits_{\Oc}a(\xi)\overline{b(\xi)}\de\xi$. 
\end{enumerate}
\end{enumerate}
\end{theorem}

\begin{proof}
See Th.~4.1.4 and Th.~2.2.7 in \cite{Pe94}. 
\end{proof}

\begin{definition}\label{moyal}
\normalfont 
Recall from Remark~\ref{smooth_oper_rem} that $\Bc(\Hc)_\infty$ is an involutive associative subalgebra of $\Bc(\Hc)$. 
It then follows by Theorem~\ref{pedersen}\eqref{pedersen_item1} 
that there exists an uniquely defined bilinear associative \emph{Moyal product} 
$$\Sc(\Oc)\times\Sc(\Oc)\to\Sc(\Oc),\quad (a,b)\mapsto a\#^\pi b $$
such that 
$$(\forall a,b\in\Sc(\Oc))\quad \Op^\pi(a\#^\pi b)=\Op^\pi(a)\Op^\pi(b). $$
Thus $\Sc(\Oc)$ is made into an involutive associative algebra 
such that the mapping $\Sc(\Oc)\to\Bc(\Hc)_\infty$, $a\mapsto\Op^\pi(a)$ 
is an algebra isomorphism. 
\qed
\end{definition}

\begin{notation}\label{distrib_vect}
\normalfont 
Recall that $\Hc_{-\infty}$ is the space of continuous antilinear functionals on $\Hc_\infty$ 
and the corresponding pairing will be denoted by 
$(\cdot\mid\cdot)\colon\Hc_{-\infty}\times\Hc_\infty\to\CC$.  
just as the scalar product in $\Hc$, since they agree on $\Hc_\infty\times\Hc_{\infty}$ 
if we think of the natural inclusions $\Hc_\infty\hookrightarrow\Hc\hookrightarrow\Hc_{-\infty}$. 
(See for instance \cite{Ca76} for more details.)
\qed
\end{notation}

\begin{definition}\label{amb}
\normalfont
If $f\in\Hc_{-\infty}$ and $\phi\in\Hc_\infty$, or $f,\phi\in\Hc$, 
then we define the corresponding \emph{ambiguity function} 
$$\Ac(f,\phi)=\Ac_\phi f\colon\gg_e\to\CC,\quad (\Ac_\phi f)(X)=(f\mid\pi(\exp_G X)\phi).$$ 
For $\phi\in\Hc_{-\infty}$ and $f\in\Hc_\infty$ we also define 
$(\Ac_\phi f)(X)=\overline{(\phi\mid\pi(\exp_G(-X))f)}$ whenever $X\in\gg_e$. 

It follows by Proposition~\ref{orthog}\eqref{orthog_item1} below that if $f,\phi\in\Hc$, 
then $\Ac_\phi f\in L^2(\gg_e)$, 
so we can use the aforementioned Fourier transform to define 
the corresponding \emph{cross-Wigner distribution} $\Wig(f,\phi)\in L^2(\Oc)$ 
such that 
$\widehat{\Wig(f,\phi)}:=\Ac_\phi f$.
\qed
\end{definition}

The second equality in Proposition~\ref{orthog}\eqref{orthog_item1} below 
could be referred to as the \emph{Moyal identity} since  
that classical identity (see for instance \cite{Gr01}) 
is recovered 
in the special case when $G$ is a simply connected Heisenberg group.

\begin{proposition}\label{orthog}
 The following assertions hold: 
\begin{enumerate}
 \item\label{orthog_item1} 
If $\phi\in\Hc$, then $\Ac_\phi f\in L^2(\gg_e)$. 
We have 
\begin{equation}\label{orthog_eq1}
\begin{aligned}
({\Ac}_{\phi_1}f_1\mid {\Ac}_{\phi_2}f_2)_{L^2(\gg_e)} 
&=(f_1\mid f_2)_{\Hc}\cdot(\phi_2\mid\phi_1)_{\Hc} \\
&=({\Wig}(f_1,\phi_1)\mid {\Wig}(f_2,\phi_2))_{L^2(\Oc)} 
\end{aligned}
\end{equation}
for arbitrary $\phi_1,\phi_2,f_1,f_2\in\Hc$. 
\item\label{orthog_item2} 
If $\phi_0\in\Hc$ with $\Vert\phi_0\Vert=1$, 
then the operator ${\Ac}_{\phi_0}\colon\Hc\to L^2(\gg_e)$, $f\mapsto {\Ac}_{\phi_0} f$, 
is an isometry and we have 
\begin{equation*}
\int\limits_{\gg_e}({\Ac}_{\phi_0}f)(X)\cdot\pi(\exp_G X)\phi\,\de X
=(\phi\mid\phi_0)f
\end{equation*}
for every $\phi\in\Hc_\infty$ and $f\in\Hc$. 
In particular, 
\begin{equation*}
\int\limits_{\gg_e}({\Ac}_{\phi_0} f)(X)\cdot\pi(\exp_G X)\phi_0\,\de X=f
\end{equation*}
for arbitrary $f\in\Hc$. 
\end{enumerate}
\end{proposition}

\begin{proof} 
See \cite{BB09c}. 
\end{proof}

\begin{corollary}\label{pseudo_prop}
The following assertions hold: 
\begin{enumerate}
\item\label{pseudo_prop_item2} 
For each $a\in\Sc(\Oc)$ we have 
$$
(\Op^\pi(a)\phi\mid f)_{\Hc}=(\widehat{a}\mid\Ac_\phi f)_{L^2(\gg_e)}
=(a\mid\Wig(f,\phi))_{L^2(\Oc)}$$
whenever $\phi,f\in\Hc$. 
Similar equalities hold if $a\in\Sc'(\Oc)$ and $\phi,f\in\Hc_\infty$. 
\item\label{pseudo_prop_item3} 
If $\phi_1,\phi_2\in\Hc_\infty$ and $a:=\Wig(\phi_1,\phi_2)\in\Sc(\Oc)$, 
then $\Op^\pi(a)$ is a rank-one operator, namely 
$\Op^\pi(a)=(\cdot \mid\phi_2) \phi_1$.  
\end{enumerate}
\end{corollary}

\begin{proof}
See \cite{BB09c}. 
\end{proof}

Assertion~\eqref{orthog_distrib_item3} in the following corollary  
in the special case of square-integrable representations 
reduces to a theorem of \cite{Co84} and \cite{CM96}. 
One thus recovers Th.~2.3 in \cite{GZ01} 
in the case of the Schr\"odinger representation 
of the Heisenberg group. 

\begin{corollary}\label{orthog_distrib}
If $\phi_0\in\Hc_\infty$ with $\Vert\phi_0\Vert=1$, 
then the following assertions hold: 
\begin{enumerate}
\item\label{orthog_distrib_item1} 
For every $f\in\Hc_{-\infty}$ we have 
\begin{equation}\label{orthog_distrib_eq1}
\int\limits_{\gg_e}({\Ac}_{\phi_0} f)(X)\cdot\pi(\exp_G X)\phi_0\,\de X=f
\end{equation} 
where the integral is convergent in the weak$^*$-topology of $\Hc_{-\infty}$. 
\item\label{orthog_distrib_item2} 
If $f\in\Hc_\infty$, then the above integral converges in 
the Fr\'echet topology of $\Hc_\infty$. 
\item\label{orthog_distrib_item3} 
If $f\in\Hc_{-\infty}$, then we have 
$f\in\Hc_\infty$ if and only if $\Ac_{\phi_0}f\in\Sc(\gg_e)$. 
\end{enumerate}
\end{corollary}

\begin{proof}
See \cite{BB09c}.
\end{proof}

\begin{remark}\label{calc_dual}
\normalfont 
Let $\Bc(\Hc)_\infty^*$ be the topological dual of the Fr\'echet space $\Bc(\Hc)_\infty$ 
and denote by $\hake{\cdot,\cdot}$ either of the duality pairings 
$$\Bc(\Hc)_\infty^*\times\Bc(\Hc)_\infty\to\CC 
\text{ and }\Sc'(\Oc)\times\Sc(\Oc)\to\CC.$$ 
Then for every tempered distribution $a\in\Sc'(\Oc)$ we can use Theorem~\ref{pedersen}\eqref{pedersen_item1} 
to define 
$\Op^\pi(a)\in\Bc(\Hc)_\infty^*$ such that  
$$(\forall b\in\Sc(\Oc))\quad \hake{\Op^\pi(a),\Op^\pi(b)}=\hake{a,b} $$
Just as in Definition~\ref{calc_def} we call $\Op^\pi(a)$ the \emph{pseudo-differential operator} with 
the \emph{symbol} $a\in\Sc'(\Oc)$. 
Note that if actually $a\in\Sc(\Oc)$, 
then the present notation agrees with Definition~\ref{calc_def} 
because of Theorem~\ref{pedersen}\eqref{pedersen_item2c}. 

The continuity properties of the above pseudo-differential operators 
can be investigated by using modulation spaces of symbols; see \cite{BB09c} for details. 
Specifically, one can introduce modulation spaces $M^{r,s}_\phi(\pi)$ 
for every unitary irreducible representation 
$\pi\colon G\to\Bc(\Hc)$.  
We always have $\Hc_\infty\subseteq M^{r,s}_\phi(\pi)$ 
and $M^{2,2}_\phi(\pi)=\Hc$. 
There exists a natural representation 
$\pi^{\#}\colon G\ltimes G\to\Bc(L^2(\Oc))$ 
such that for suitable $\Phi\in\Sc(\Oc)\setminus\{0\}$, 
the Weyl calculus $\Op^\pi(\cdot)$ defines a continuous linear mapping 
from the modulation space $M^{\infty,1}_\Phi(\pi^{\#})$ into the space of \emph{bounded} linear operators 
on~$\Hc$. 
One of the main theorems of \cite{GH99} is recovered in the special case when $\pi$ is the Schr\"odinger representation 
of the $(2n+1)$-dimensional Heisenberg group. 
Some new results related to this circle of ideas 
will be established in Section~\ref{Sect5} below. 
\qed 
\end{remark}

\begin{remark}[\cite{BB09c}]\label{wig_ext}
\normalfont
We can define the cross-Wigner distribution $\Wig(f_1,f_2)\in\Sc'(\Oc)$ 
for arbitrary $f_1,f_2\in\Hc_{-\infty}$ as follows. 
By using Corollary~\ref{regularization} we can define for $f_1,f_2\in\Hc_{-\infty}$
the continuous antilinear functional
$$T_{f_1,f_2}\colon\Bc(\Hc)_\infty\to\CC,\quad 
T_{f_1,f_2}(A):=(f_1\mid Af_2).$$
That is, $T_{f_1,f_2}\in\Bc(\Hc)_\infty^*$, 
and then Th.~4.1.4(5) in \cite{Pe94} shows that 
there exists a unique distribution $a_{f_1,f_2}\in\Sc'(\Oc)$ 
such that $\Op^\pi(a_{f_1,f_2})=T_{f_1,f_2}$. 
Now define 
$$\Wig(f_1,f_2):=a_{f_1,f_2}.$$ 
We can consider the rank-one operator 
$S_{f_1,f_2}:=(\cdot\mid f_2)f_1\colon \Hc_\infty\to\Hc_{-\infty}$ and 
for arbitrary $A\in\Bc(\Hc)_\infty$ thought of as 
a continuous linear map $A\colon\Hc_{-\infty}\to\Hc_\infty$ as above 
we have 
$$\Tr(S_{f_1,f_2}A)=(f_1\mid Af_2)=T_{f_1,f_2}(A).$$ 
Thus the trace duality pairing allows us to identify 
the functional $T_{f_1,f_2}\in\Bc(\Hc)_\infty^*$ 
with the rank-one operator $(\cdot\mid f_2)f_1$, 
and then we can write 
\begin{equation}\label{wig_ext_eq1}
(\forall f_1,f_2\in\Hc_{-\infty})\quad 
\Op^\pi(\Wig(f_1,f_2))=(\cdot\mid f_2)f_1.
\end{equation}
In particular, it follows that the above extension of 
the cross-Wigner distribution to a mapping 
$\Wig(\cdot,\cdot)\colon\Hc_{-\infty}\times\Hc_{-\infty}\to\Sc'(\Oc)$ 
allows us to generalize 
the assertion of Corollary~\ref{pseudo_prop}\eqref{pseudo_prop_item3} 
to arbitrary $\phi_1,\phi_2\in\Hc_{-\infty}$. 
\qed
\end{remark}

\section{Modulation spaces}\label{Sect5}

The modulation spaces play a central role in the time-frequency analysis 
(see \cite{Gr01}) and proved to be a very useful tool in the study
of continuity properties of pseudo-differential operators (\cite{GH99}). 
These classical ideas can be formulated within the representation theory 
of the Heisenberg groups, and this representation theoretic viewpoint 
turned out to be very effective in order to extend the corresponding notions 
to unitary irreducible representations of arbitrary nilpotent Lie groups 
(see \cite{BB09c}). 
In the first two subsections of the present section we shall provide 
some preparations and then describe the general notion of 
modulation spaces introduced in \cite{BB09c}. 
We eventually illustrate this notion by discussing a specific class 
of irreducible representations on Hilbert spaces of the form $L^2(\Oc)$, 
where $\Oc$ is any coadjoint orbit of a nilpotent Lie group 
(see Proposition~\ref{double} and Remark~\ref{double_cor2}). 

\subsection{Semidirect products}

\begin{definition}\label{sd_def}
\normalfont
Let $G_1$ and $G_2$ be connected Lie groups and assume that we have 
 a continuous group homomorphism 
$\alpha\colon G_1\to\Aut G_2$, $g_1\mapsto\alpha_{g_1}$. 
The corresponding \emph{semidirect product of Lie groups} 
$G_1\ltimes_\alpha G_2$ is the connected Lie group whose 
underlying manifold is the Cartesian product $G_1\times G_2$ 
and whose group operation is given by 
\begin{equation}\label{sd_prod}
(g_1,g_2)\cdot(h_1,h_2)=(g_1h_1,\theta_{h_1^{-1}}(g_2)h_2) 
\end{equation}
whenever $g_j,h_j\in G_j$ for $j=1,2$. 

Let us denote by $\dot\alpha\colon\gg_1\to\Der\gg_2$ 
the homomorphism of Lie algebras defined as the differential 
of the Lie group homomorphism $G_1\to\Aut\gg_2$, $g_1\mapsto\Lie(\alpha_{g_1})$. 
Then the \emph{semidirect product of Lie algebras} 
$\gg_1\ltimes_{\dot\alpha}\gg_2$ 
is the Lie algebra whose underlying linear space is the Cartesian product 
$\gg_1\times\gg_2$ with the Lie bracket given by 
\begin{equation}\label{sd_bracket}
[(X_1,X_2),(Y_1,Y_2)]=([X_1,Y_1],\dot\alpha(X_1)Y_2-\dot\alpha(Y_1)X_2+[X_2,Y_2])
\end{equation}
if $X_j,Y_j\in\gg_j$ for $j=1,2$. 
One can prove that $\gg_1\ltimes_{\dot\alpha}\gg_2$ is the Lie algebra 
of the Lie group $G_1\ltimes_\alpha G_2$ 
(see for instance Ch.~9 in \cite{Ho65}). 
\qed
\end{definition}

\begin{remark}\label{sd_rem}
\normalfont
Assume that $G_1$ and $G_2$ are nilpotent Lie groups and  
$\alpha\colon G_1\to\Aut G_2$ is a \emph{unipotent automorphism}.  
That is, for every $X_1\in\gg_1$ there exists an integer $m\ge1$ such that 
$\dot\alpha(X_1)^m=0$. 
Then an inspection of \eqref{sd_bracket} shows that $\gg_1\ltimes_{\dot\alpha}\gg_2$ 
is a nilpotent Lie algebra, hence $G_1\ltimes_\alpha G_2$ is a nilpotent Lie group. 
\qed
\end{remark}

\begin{example}\label{sd_ex1}
\normalfont
Let $G$ be a nilpotent Lie group. 
If we specialize Definition~\ref{sd_def} for $G_1:=G$, $G_2=(\gg,+)$, and 
$\alpha:=\Ad_G\colon G\to\Aut\gg$, then we get 
the semidirect product $G\ltimes_{\Ad_G}\gg$ which is a nilpotent Lie group by 
Remark~\ref{sd_rem} and is isomorphic to the tangent group $TG$. 
The Lie algebra of $G\ltimes_{\Ad_G}\gg$ is $\gg\ltimes_{\ad_{\gg}}\gg_0$  
(where $\gg_0$ stands for the abelian Lie algebra that has the same 
underlying linear space as~$\gg$)
and the corresponding exponential map is given by 
$$\begin{aligned}
\exp_{G\ltimes_{\Ad_G}\gg}(X,Y)
&=(\exp_G X,\int\limits_0^1\Ad_G(\exp_G(sX))Y\,\de s) \\
&=(\exp_G X,\int\limits_0^1\ee^{s\cdot\ad_{\gg}X}Y\,\de s)
\end{aligned}$$ 
for every $(X,Y)\in\gg\ltimes_{\ad_{\gg}}\gg_0$ 
(see for instance Prop.~2.7(2) in \cite{BB09a}). 
\qed
\end{example}

\subsection{Modulation spaces for unitary representations}
In this short subsection we just recall the definition 
of the modulation spaces for the unitary irreducible representations of nilpotent Lie groups. 
We refer to \cite{BB09c} for a more detailed discussion of this notion. 

\begin{definition}\label{modular_def}
\normalfont
Let $\phi\in\Hc_\infty\setminus\{0\}$ be fixed and 
assume that we have a direct sum decomposition 
$\gg_e=\gg_e^1\dotplus\gg_e^2$.

Then let $1\le r,s\le\infty$ and   
for arbitrary $f\in\Hc_{-\infty}$ define 
$$\Vert f\Vert_{M^{r,s}_\phi}
=\Bigl(\int\limits_{\gg_e^2}
\Bigl(\int\limits_{\gg_e^1}
\vert(\Ac_\phi f)(X_1,X_2)\vert^r\de X_1 \Bigr)^{s/r}
\de X_2\Bigr)^{1/s}\in[0,\infty] $$
with the usual conventions if $r$ or $s$ is infinite. 
Then we call the space 
$$M^{r,s}_\phi(\pi):=\{f\in\Hc_{-\infty}\mid\Vert f\Vert_{M^{r,s}_\phi}<\infty\}$$ 
a \emph{modulation space} for the irreducible unitary representation $\pi\colon G\to\Bc(\Hc)$ 
with respect to the decomposition $\gg_e\simeq\gg_e^1\times\gg_e^2$ 
and the \emph{window vector} $\phi\in\Hc_\infty\setminus\{0\}$. 
\qed
\end{definition}

\begin{example}\label{mod_L2}
\normalfont
For any choice of $\phi\in\Hc_\infty\setminus\{0\}$ in Definition~\ref{modular_def} we have 
$$M^{2,2}_\phi(\pi)=\Hc.$$
Indeed, this equality holds since 
$\Vert\Ac_\phi f\Vert_{L^2(\gg_e)}=\Vert\phi\Vert\cdot\Vert f\Vert$ 
for every $f\in\Hc$ 
(see Proposition~\ref{orthog} above). 
\qed
\end{example}

\subsection{A specific irreducible representation on $L^2(\Oc)$} 
We are going to construct here some irreducible representations on the Hilbert spaces of the form $L^2(\Oc)$, 
where $\Oc$ can be any coadjoint orbit of a nilpotent Lie group. 
A different construction involving the Moyal product (see Definition~\ref{moyal}) 
was used in Def.~2.19 in the paper \cite{BB09c} in order to get a representation $\pi^\#$ 
with the same representation space $L^2(\Oc)$. 
The modulation spaces for $\pi^\#$ 
turned out to be relevant for establishing 
the continuity properties of the pseudo-differential operators 
obtained by the Weyl-Pedersen calculus for any unitary representation 
associated with the coadjoint orbit~$\Oc$.
(See also Remark~\ref{calc_dual}.)

\begin{proposition}\label{double}
Let $Z$ be the center of the connected, simply connected, nilpotent Lie group~$G$ 
with the corresponding Lie algebra $\zg\subseteq\gg$. 
Endow the coadjoint orbit~$\Oc$ with a Liouville measure and define 
$$\widetilde{\pi}\colon G\ltimes_{\Ad}\gg\to\Bc(L^2(\Oc)),\quad 
(\widetilde{\pi}(g,Y)f)(\xi)=\ee^{\ie \hake{\xi,Y}}f(\Ad^*_G(g^{-1})\xi).$$
Then the following assertions hold: 
\begin{enumerate}
\item\label{double_item1} 
The group $\widetilde{G}:=G\ltimes_{\Ad}\gg$ 
is nilpotent and its center is $Z\times\zg$. 
\item\label{double_item2} 
$\widetilde{\pi}$ is a unitary irreducible representation 
of $\widetilde{G}$. 
\item\label{double_item2and1/2} 
Let us denote by $\widetilde{\gg}=\gg\ltimes_{\ad_{\gg}}\gg_0$ 
the Lie algebra of $\widetilde{G}$ 
(where $\gg_0$ stands for the abelian Lie algebra with the same 
underlying linear space as~$\gg$)
and define 
$$\widetilde{X}_j=
\begin{cases}
\hfill (0,X_j) &\text{ for }j=1,\dots,n,\\
(X_{j-n},0) &\text{ for }j=n+1,\dots,2n.
\end{cases} 
$$
Then $\widetilde{X}_1,\dots,\widetilde{X}_{2n}$ is a Jordan-H\"older basis 
in $\widetilde{\gg}$ and the corresponding predual for 
the coadjoint orbit $\widetilde{\Oc}\subseteq\widetilde{\gg}^*$ associated with 
the representation $\widetilde{\pi}$ is $$\widetilde{\gg}_{\widetilde{e}}=\gg_e\times\gg_e\subseteq\widetilde{\gg},$$
where $\widetilde{e}$ is the set of jump indices for~$\widetilde{\Oc}$. 
\item\label{double_item4} 
The space of smooth vectors for the representation $\widetilde{\pi}$ 
is $\Sc(\Oc)$. 
\end{enumerate}
\end{proposition}

\begin{proof}
\eqref{double_item1}
Recall that the multiplication in the semi-direct product group $\widetilde{G}$ 
is given by 
$$(g_1,Y_1)\cdot(g_2,Y_2)=(g_1g_2,Y_1+\Ad_G(g_1)Y_2) $$
while the bracket in the corresponding Lie algebra $\widetilde{\gg}=\gg\ltimes_{\ad}\gg$ is 
defined by 
\begin{equation}\label{double_eq0}
[(X_1,Y_1),(X_2,Y_2)]=([X_1,X_2],[X_1,Y_2]-[X_2,Y_1]).
\end{equation}
An inspection of these equations quickly leads to the conclusion that $\widetilde{\gg}$ 
is a nilpotent Lie algebra with the center $\zg\times\zg$. 

\eqref{double_item2} 
If $(g_1,Y_1),(g_2,Y_2)\in\widetilde{G}$ and $f\in L^2(\Oc)$, 
then for $\xi\in\Oc$ we have 
$$\begin{aligned}
\widetilde{\pi}(g_1,Y_1)(\widetilde{\pi}(g_2,Y_2)f)(\xi)
&=\ee^{\ie\hake{\xi,Y_1}}(\widetilde{\pi}(g_2,Y_2)f)(\Ad_G^*(g_1^{-1})\xi) \\
&=\ee^{\ie\hake{\xi,Y_1}}\ee^{\ie\hake{\Ad_G^*(g_1^{-1})\xi,Y_2}}f(\Ad_G^*(g_2^{-1})\Ad_G^*(g_1^{-1})\xi) \\
&=\ee^{\ie\hake{\xi,Y_1+\Ad_G(g_1)Y_2}}f(\Ad_G^*((g_1g_2)^{-1})\xi) 
\end{aligned}$$
hence $\widetilde{\pi}(g_1,Y_1)\widetilde{\pi}(g_2,Y_2)=\widetilde{\pi}((g_1,Y_1)(g_2,Y_2))$. 
Next note that the representation $\widetilde{\pi}$ is unitary 
since the Liouville measure on $\Oc$ is invariant under the coadjoint 
action of $G$. 

To see that $\widetilde{\pi}$ is irreducible, let $T\colon L^2(\Oc)\to L^2(\Oc)$ 
be any bounded linear operator satisfying $T\widetilde{\pi}(g,Y)=\widetilde{\pi}(g,Y)T$ 
for arbitrary $(g,Y)\in\widetilde{G}$. 
We have to check that $T$ is a scalar multiple of the identity operator on $L^2(\Oc)$. 
By applying the assumption for $g=\1\in G$ we see that 
$T$ belongs to the commutant of the family of multiplication operators 
by the functions in the family $\{\ee^{\ie\hake{Y,\cdot}}\mid Y\in\gg\}\subseteq L^\infty(\Oc)$. 
On the other hand, 
we recall that the mapping 
$$\Oc\to\RR^d,\quad \xi\to(\hake{\xi,X_{j_1}},\dots,\hake{\xi,X_{j_d}}) $$
is a global chart which takes the Liouville measure of $\Oc$ 
to a Lebesgue measure on~$\RR^d$
(see for Remark~\ref{jump_remark}). 
Then we can use the Fourier transform to see that the linear subspace generated by $\{\ee^{\ie\hake{Y,\cdot}}\mid Y\in\gg\}$ 
is weak$^*$-dense in $L^\infty(\Oc)$ ($\simeq L^1(\Oc)^*$). 
Therefore the operator $T\colon L^2(\Oc)\to L^2(\Oc)$ commutes with 
all the multiplication operators by functions in $L^\infty(\Oc)$, 
and then it has to be in turn the multiplication operator by some function $\phi\in L^\infty(\Oc)$. 
Now, by using the assumption that $T$ commutes with $\pi(g,0)$ for arbitrary $g\in G$,  
it easily follows that $\phi$ has to be a constant function 
since the coadjoint action of $G$ on the orbit $\Oc$ is transitive.

\eqref{double_item2and1/2} 
It is straightforward to check that 
$\widetilde{X}_1,\dots,\widetilde{X}_{2n}$ is a Jordan-H\"older basis 
in $\widetilde{\gg}$. 
Next note that $\Sc(\Oc)$ is contained in the space of smooth vectors 
for the representation $\widetilde{\pi}$ and for arbitrary 
$f\in\Sc(\Oc)$ and $(X,Y)\in\widetilde{\gg}$ 
we have 
\begin{equation}\label{double_proof_eq1}
(\forall\xi\in\Oc)\quad 
(\de\widetilde{\pi}(X,Y)f)(\xi)=\ie\hake{\xi,Y}f(\xi)
+\frac{\de}{\de t}\Big{\vert}_{t=0}f(\xi\circ\ee^{t\ad_{\gg}X}). 
\end{equation}
It then follows by a straightforward application of Prop.~2.4.1 in \cite{Pe84} 
and by Lemmas 1.4.1 and 1.5.1 in \cite{Pe89} 
that the set of jump indices for the coadjoint orbit $\widetilde{\Oc}$ 
is $\widetilde{e}=\{j_1,\dots,j_d,n+j_1,\dots,n+j_d\}$, and then 
$\widetilde{\gg}_{\widetilde{e}}=\gg_e\times\gg_e\subseteq\widetilde{\gg}$. 

\eqref{double_item4} 
It follows by \eqref{double_proof_eq1} and 
by Lemmas 1.4.1 and 1.5.1 in \cite{Pe89} 
again that there exists a polynomial chart on $\Oc$ such that 
in the corresponding chart, the associative algebra 
generated by the image of $\de\widetilde{\pi}$ 
contains all the linear partial differential operators with polynomial coefficients. 
This implies that the space of smooth vectors for the representation $\widetilde{\pi}$ 
is equal to $\Sc(\Oc)$, as asserted.
\end{proof}

\begin{corollary}\label{double_cor}
Assume the setting of Proposition~\ref{double}. 
The ambiguity function 
$$\widetilde{\Ac}\colon L^2(\Oc)\times L^2(\Oc)\to L^2(\widetilde{\gg}_{\widetilde{e}})
=L^2(\gg_e\times\gg_e) $$
for the representation $\widetilde{\pi}\colon\widetilde{G}\to\Bc(L^2(\Hc))$ 
is given by the formula 
$$(\widetilde{\Ac}_\phi f)(X,Y)
=\int\limits_{\Oc}\ee^{-\ie\hake{\xi,\int\limits_0^1\ee^{s\cdot\ad_{\gg}X}Y\,\de s}} 
f(\xi)\overline{\phi(\xi\circ\ee^{\ad_{\gg}X})} \de\xi$$ 
for arbitrary $X,Y\in\gg_e$ and $f,\phi\in L^2(\Oc)$. 
\end{corollary}

\begin{proof}
For every $f,h\in L^2(\Oc)$ we have 
\begin{equation}\label{double_cor_proof_eq1}
(\widetilde{\Ac}_\phi f)(X,Y)=(f\mid\widetilde{\pi}(\exp_{\widetilde{G}}(X,Y))\phi)_{L^2(\Oc)}.
\end{equation}
On the other hand, for the element $(X,Y)\in\widetilde{\gg}$ we have 
$$\exp_{\widetilde{G}}(X,Y)=(\exp_G X,\int\limits_0^1\Ad_G(\exp_G(sX))Y\,\de s)
=(\exp_G X,\int\limits_0^1\ee^{s\cdot\ad_{\gg}X}Y\,\de s)$$
(see Example~\ref{sd_ex1})
hence 
$$\begin{aligned}
(\widetilde{\pi}(\exp_{\widetilde{G}}(X,Y))\phi)(\xi)
&=\ee^{\ie\hake{\xi,\int\limits_0^1\ee^{s\cdot\ad_{\gg}X}Y\,\de s}} 
\phi(\Ad_G^*(\exp_G(-X))\xi) \\
&=\ee^{\ie\hake{\xi,\int\limits_0^1\ee^{s\cdot\ad_{\gg}X}Y\,\de s}} 
\phi(\xi\circ\ee^{\ad_{\gg}X})
\end{aligned}$$
and then the conclusion follows by \eqref{double_cor_proof_eq1}. 
\end{proof}

\begin{remark}\label{double_cor2}
\normalfont
Assume the setting of the above Proposition~\ref{double}. 
It follows by Corollary~\ref{double_cor} along with 
Schur's criterion for integral operators that there exists 
a constant $C_\Phi>0$ such that for every 
$F\in L^2(\Oc)$ and $Y\in\gg_e$ we have 
$\Vert\widetilde{\Ac}_\Phi F(\cdot,Y)\Vert_{L^2(\gg_e)}
\le C_\Phi\Vert F\Vert_{L^2(\Oc)}$, hence 
$\Vert F\Vert_{M^{2,\infty}_\Phi(\widetilde{\pi})}\le 
C_\Phi \Vert F\Vert_{L^2(\Oc)}$. 
Therefore there exists a continuous inclusion map 
$L^2(\Oc)\hookrightarrow M^{2,\infty}_\Phi(\widetilde{\pi})$. 
See also \cite{BB09b} for similar inclusion maps for the modulation spaces 
in the setting of the magnetic Weyl calculus on nilpotent Lie groups.
\qed
\end{remark}

\begin{example}\label{two-step}
\normalfont 
Assume that $\gg$ is two-step nilpotent Lie algebra, that is, we have $[\gg,[\gg,\gg]]=\{0\}$. 
Let $\Oc\subseteq\gg^*$ be any nontrivial coadjoint orbit 
and pick $\xi_0\in\Oc$. 
If we denote by $\zg$ the center of $\gg$, then 
$$\Oc=\{\xi\in\gg^*\mid\quad \xi\vert_{\zg}=\xi_0\vert_{\zg}\},$$
since $\Oc$ is a flat orbit. 
Then by Corollary~\ref{double_cor} along with the fact that $[\gg,\gg]\subseteq\zg$ we get 
$$\begin{aligned}
(\widetilde{\Ac}_\phi f)(X,Y)
&=\int\limits_{\Oc}\ee^{-\ie\hake{\xi,\int\limits_0^1\ee^{s\cdot\ad_{\gg}X}Y\,\de s}} 
f(\xi)\overline{\phi(\xi\circ\ee^{\ad_{\gg}X})} \de\xi \\
&=\int\limits_{\Oc}\ee^{-\ie\hake{\xi,Y+\frac{1}{2}[X,Y]}} 
f(\xi)\overline{\phi(\xi+\xi\circ\ad_{\gg}X)} \de\xi \\
&=\ee^{-\frac{\ie}{2}\hake{\xi_0,[X,Y]}}\int\limits_{\Oc}\ee^{-\ie\hake{\xi,Y}} 
f(\xi)\overline{\phi(\xi+\xi_0\circ\ad_{\gg}X)} \de\xi. 
\end{aligned}$$ 
The above formula and suitable global coordinates on $\Oc$, 
one shows that the ambiguity function of the representation $\widetilde{\pi}$ 
agrees with the ambiguity function of the Schr\"odinger representation of a certain Heisenberg group, 
as defined in \cite{Gr01}. 
\qed
\end{example}

\textbf{Acknowledgment.} 
Partial financial support from the grant PNII - Programme ``Idei'' (code 1194) 
is acknowledged.

\end{document}